\documentclass[final]{amsart}
\usepackage{pdfsync}
\usepackage{amsmath,amssymb}
\usepackage{enumerate}
\usepackage{xspace}
\usepackage{boxedminipage}
\usepackage{algorithmicx}
\usepackage[ruled]{algorithm}
\usepackage{algpseudocode}
\usepackage{listings}
\usepackage{tabularx}
\usepackage{subfigure}
\usepackage{tristan}
\usepackage{fullpage}

\renewcommand{\vec}[1]{\geovec{#1}}
\renewcommand{\mat}[1]{\geomat{#1}}

\renewcommand{\enorm}[2]{\ensuremath{\Norm{#1}}_{dG,{#2}}}

\renewcommand{\leb}[1]{\ensuremath{\LL^{#1}}}

\newcommand{\figdir}{lowresfigures/}

\numberwithin{equation}{section}
\setboolean{showtodo}{true}
\setboolean{shownotes}{true}
\setboolean{showchanges}{false}
\setboolean{usemathrsfs}{false}
\usepackage{tkz-graph}
\usetikzlibrary{arrows}

\author{
  James Jackaman
}
\address{
  James Jackaman
  \thanks{
    Department of Mathematics and Statistics, University of Reading, Reading RG6 6AX, UK.
    {\tt{James.Jackaman@pgr.reading.ac.uk}}.
}}

\author{
  Georgios Papamikos
}
\address{
  Georgios Papamikos
  \thanks{
    Department of Mathematics and Statistics, University of Reading, Reading RG6 6AX, UK.
    {\tt{G.Papamikos@reading.ac.uk}}.
}}

\author{
  Tristan Pryer
}
\address{
  Tristan Pryer
  \thanks{
    Department of Mathematics and Statistics, University of Reading, Reading RG6 6AX, UK.
    {\tt{T.Pryer@reading.ac.uk}}.
}}

\thanks{JJ was supported through a PhD scholarship awarded by the ``EPSRC Centre for Doctoral Training in the Mathematics of Planet Earth at Imperial College London and the University of Reading'' EP/L016613/1.
  TP was partially supported through the EPSRC grant ``Adaptive regularisation'' EP/P000835/1.
  GP was supported through a Rado Postdoctoral fellowship at the University of Reading. All this support is gratefully acknowledged.}

\title[Conservative discretisations for the vectorial modified KdV system]
      {The design of conservative finite element discretisations for the vectorial modified KdV equation}
\date{\today}

\begin{document}
\maketitle
\begin{abstract}
  We design a consistent Galerkin scheme for the approximation of the
  vectorial modified Korteweg-de Vries equation. We demonstrate that
  the scheme conserves energy up to machine precision. In this sense
  the method is consistent with the energy balance of the continuous
  system. This energy balance ensures there is no numerical
  dissipation allowing for extremely accurate long time simulations
  free from numerical artifacts. Various numerical experiments are
  shown demonstrating the asymptotic convergence of the method with
  respect to the discretisation parameters. Some simulations are also
  presented that correctly capture the unusual interactions between
  solitons in the vectorial setting.
\end{abstract}

\section{Introduction}
\label{sec:introduction}

Hamiltonian partial differential equations (PDEs) are a specific class
of PDE endowed with physically relevant algebraic and geometric
structures \cite{Olver:1993}. They arise form a variety of areas
\cite{Shepherd:1990}, not least meteorological, such as the
semi-geostrophic equations \cite{RoulstoneNorbury:1994}, and
oceanographical, such as the Korteweg-de Vries (KdV) and nonlinear
Schr\"odinger equations \cite{MullerGarretOsborne:2005}. The KdV and
nonlinear Schr\"odinger equations are particularly special examples,
in that they are bi-Hamiltonian \cite{Magri:1978}. This means they
have two different Hamiltonian formulations which, in turn, is one way
to understand the notion of \emph{integrability} of these
problems. Regardless, the applications and the need to quantify the
dynamics of the general Cauchy problem motivate the development of
accurate long time simulations for reliable prediction of dynamics in
both meteorology and oceanography.

A difficulty in the design of schemes for this class of problems is
that the long term dynamics of solutions can be destroyed by the
addition of \emph{artificial numerical diffusion}. The reason for
inclusion of this in a given scheme is the desirable stability
properties this endows on the approximation, however, this typically
destroys all information in the long term dynamics of the system
through smearing of solutions.

Conservative schemes for Hamiltonian ordinary differential equations
(ODEs) are relatively well understood, see
\cite[c.f.]{Celledoni:2012,LeimkuhlerReich:2004,HairerLubichWanner:2006,BokhoveLynch:2007,BlanesCasa:2016}. Typically
numerical schemes designed for this class of problem which have some
property of the ODE built into them, for example conservativity of the
Hamiltonian or the underlying symplectic form, are classified as
geometric integrators. In the PDE setting symplectic structure
preserving algorithms have also been developed
\cite{Reich:1999,BridgesReich:2001,CohenOwrenRaynaud:2008,CohenMatsuoRaynaud:2014}. Here
the PDEs are rewritten using their corresponding symplectic form, and
the notion of structure preservation is given in terms of a discrete
(difference) conservation law involving differential forms. Typically
these problems are solved using an Euler/Preissman box scheme, although
space-time finite element methods under an appropriate quadrature
choice form a natural generalisation.


In this contribution we consider a system of equations that are a
multidimensional generalisation of the famous modified KdV equation
\begin{equation*}
  u_t
  +
  \frac{3}{2}
  u^2
  u_x
  +
  u_{xxx} = 0,
  \label{eq:mkdv}
\end{equation*}
where the subindices denote partial differentiation with respect to
the corresponding independent variable. This equation has numerous
applications not least including fluid dynamics and plasma physics
\cite{AblowitzClarkson:1991}. The vectorial equation
\cite{Mari-BeffaSandersWang:2002} appears frequently in the study of
ocean waves and Riemannian geometry
\cite{SandersWang:2003,AncoNgatatWilloughby:2011,Anco:2006}. The
presence of the Korteweg third order term as well as the lack of sign
in the corresponding energy functional, further details given in
\S\ref{sec:vmkdv}, can cause issues in the numerical treatment of this
problem. The vectorial equation also has additional complications that
do not arise in the scalar setting. Indeed, the underlying
conservation laws themselves have a vastly different structure, for
example there is no mass conservation in the vectorial case. This is
due to the different form of Hamiltonian operator.

In previous numerical studies of the scalar KdV and modified KdV
equations \cite[c.f.]{XuShu:2007,YanShu:2002}, it has been observed
that classical finite volume and discontinuous Galerkin (dG) schemes with
``standard'' numerical fluxes introduce numerical artifacts. These
typically appear through numerical regularisation effects included for
stability purposes that are not, however, adapted to the variational
structure of the problem. The result of these artifacts is an inconsistency in
the discrete energy.

Hamiltonian problems are inherently conservative, that is, the
underlying Hamiltonian is conserved over time. Other equations,
including those of integrable type, may have additional structures
that manifest themselves through additional conserved quantities. In
particular, mass and momentum are such quantities. In
\cite{BonaChenKarakashianXing:2013} and
\cite{KarakashianMakridakis:2015} the authors propose and analyse a dG
scheme for generalised KdV equations. The scheme itself is very
carefully designed to be conservative, in that the invariant
corresponding to the \emph{momentum} is inherited by the
discretisation. This yields $\leb{2}$ stability quite naturally in the
numerical scheme and extremely good long time dynamics. In the scalar
case one may also design schemes that conserve the energy itself
\cite{Winther:1980,JackamanPryer:2017}, however, it does not seem
possible to design schemes to conserve more than two of these
invariants. A natural question is then ``conservation of which
invariant yields the best scheme?'' \cite{Jackaman:2017}.

When examining the case of systems of Hamiltonian equations much less
work has been carried out, for example \cite{BaiZhang:2009} give a
near conservative method for a system of Schr\"odinger--KdV type and
\cite{BonaDougalisMitsotakis:2007} study a system of KdV equations. To
the authors knowledge there has not been any work on the system we
study, nor such Hamiltonian systems in general. Our goal in this work
is the derivation of Galerkin discretisations aimed at preserving the
underlying algebraic properties satisfied by the PDE system whilst
avoiding the introduction of stabilising diffusion terms. Our schemes
are therefore consistent with (one of) the Hamiltonian formulation of
the original problem. It is important to note that our approach is not
an adaptation of entropy conserving schemes developed for systems of
conservation laws, rather we study the algebraic properties of the PDE
and formulate the discretisation to inherit this specific
structure. The methods are of arbitrarily high order of accuracy in
space and provide relevant approximations free from numerical
artifacts. To the authors knowledge this is the first class of finite
element method for this particular class of system to have these
properties.

The rest of this work is set out as follows: In \S\ref{sec:vmkdv} we
introduce notation, the model problem and some of its properties. We
give a reformulation of the system through the introduction of
auxiliary variables, these are introduced to allow for a simple
construction of the numerical scheme. In \S\ref{sec:temp} we examine a
temporal discretisation of the problem that guarantees the
conservation of the underlying Hamiltonian. In
\S\ref{sec:spat} we propose a spatial discretisation based on
continuous finite elements, state a fully discrete method and show
that it is conservative. Finally in \S\ref{sec:numerics} we summarise
extensive numerical experiments aimed at testing the robustness of the
method in long time simulations.

\newcommand{\vmkdv}{vmKdV\xspace}

\section{The vectorial modified KdV equation}
\label{sec:vmkdv}

In this section we formulate the model problem, fix notation and give
some basic assumptions. We describe some known results and history of
the vectorial modified Korteweg-de Vries (\vmkdv) equation,
highlighting the Hamiltonian structure of the equation. We show that
the underlying Hamiltonian structure naturally yields an induced
stability of the solutions to the PDE system and give a brief
description of how to construct some exact solutions using a dressing
method. We then show how the system can be written through induced
auxiliary variables which are the basis of the design of our
numerical scheme.

Throughout this work we denote the standard
Lebesgue spaces by $\leb{p}(\w)$ for $\w\subseteq \reals$, $p\in
[1,\infty]$, equipped with corresponding norms
$\Norm{u}_{\leb{p}(\w)}$. In addition, we will denote $\sobh{k}(\w)$
to be the Hilbert space of order $k$ of real-valued functions defined
over $\w\subseteq \reals$ with norm $\Norm{u}_{\sobh{k}(\w)}$.

The \vmkdv equation is an evolutionary PDE for a real, $d$-vector
valued function
\begin{equation}
  \dfunkmapsto
      {\vec u}
      {(x,t)}
      {\reals^2}
      {\vec u(x,t) = \Transpose{\qp{u_1 \dots, u_d}}}
      {\reals^d}
\label{eq:udef}
\end{equation}
and is given by
\begin{equation}
  \vec u_t
  +
  \frac{3}{2}
  {\vec u}\cdot{\vec u}
  \vec u_x
  +
  \vec u_{xxx} = \vec 0.
\label{eq:vmkdv}
\end{equation}
Here we are using ``$\vec x\cdot \vec y$'' as the Euclidean inner
product between two vectors, $\vec x$ and $\vec y$ and ``$\norm{\vec
  x}$'' as the induced Euclidean norm of $\vec x$.

A particular case of the \vmkdv system occurs when $d=2$, $\vec u =
\Transpose{\qp{u_1, u_2}}$ when \eqref{eq:vmkdv} can be identified
with the complex modified KdV (mKdV) equation
\begin{equation}
  y_t+\frac{3}{2}\norm{y}^2y_x+y_{xxx}=0
  \label{eq:cmkdv}
\end{equation}
for the complex dependent variable $y=u_1+ i u_2$. Sometimes this is
also called the Hirota mKdV equation \cite{Hirota:1973}. When $d=1$ we
obtain the famous mKdV equation which has been studied numerically in
the context of Galerkin methods in
\cite{BonaChenKarakashianXing:2013,JackamanPryer:2017}.

Equation \eqref{eq:vmkdv} admits both Lie and discrete point
symmetries and has infinitely many conservation laws. Indeed, under
the action of the orthogonal group $O_d(\mathbb{R})$, which is the
group of $d\times d$ matrices such that $\Transpose{\geomat A}\geomat
A=\eye$,
\begin{equation}
\widetilde{\vec{u}}=\geomat A\vec{u}, \quad \text{ for }\geomat A\in O_d(\mathbb{R})
\label{eq:O(N)}
\end{equation}
equation \eqref{eq:vmkdv} remains invariant. Moreover, vmKdV is invariant under the translations
\begin{equation}
\widetilde{x}=x+\epsilon, \quad \widetilde{t}=t+\gamma
\label{eq:translations}
\end{equation}
and under the scaling transformation
\begin{equation}
\left(\widetilde{x},\widetilde{t},\widetilde{\vec{u}}\right)=\left(e^{\epsilon}x,e^{3\epsilon}t,e^{-\epsilon}\vec{u}\right).
\label{eq:scaling}
\end{equation}  

\begin{Pro}[Conservative properties of solutions]
  \label{the:cons}
  The \vmkdv equation admits the following conservation laws:
  \begin{equation}
    \begin{split}
      \D_t f_2(\vec u) &= \D_x g_2(\vec u)
      \\
      \D_t f_4(\vec u) &= \D_x g_4(\vec u),
    \end{split}
  \end{equation}
  where the conserved densities are given by
  \begin{equation}
    \begin{split}
      f_2(\vec u) &= \frac 1 2 \norm{\vec u}^2
      \\
      f_4(\vec u) &= \frac 1 2 \norm{\vec u_x}^2 - \frac 1 8 \norm{\vec u}^4
    \end{split}
  \end{equation}
  and the corresponding fluxes are
  \begin{equation}
    \begin{split}
      g_2(\vec u) &= \norm{\vec u_x}^2 - 2 \vec u \cdot \vec u_{xx} - \frac 3 4 \norm{\vec u}^4
      \\
      g_4(\vec u) &= \frac 1 8 \norm{\vec u}^6 - \norm{\vec u}^2 \norm{\vec u_x}^2 - \frac 1 2 \qp{\vec u \cdot \vec u_x}^2 - \vec u_x \cdot \vec u_{xxx} + \frac 1 2 \norm{\vec u_{xx}}^2 + \frac 1 2 \norm{\vec u}^2 \vec u\cdot \vec u_{xx}.
    \end{split}
  \end{equation}
  \end{Pro}
\begin{Proof}
To prove that the total time derivative of $f_2(\textbf{u})$ and $f_4(\textbf{u})$ are in the image of $D_x$ we use the Euler operator $\vec{E}=(E_1,...,E_d)$, where
$$
E_i(f)=\sum_{k=0}^{\infty}(-D_x)^k \partial_{u_{i_{kx}}}f, \quad u_{i_{kx}}=u_{i_{\underbrace{x...x}_k}}
$$
and the fact that $\text{Ker}\vec E=\text{Im}D_x$, see \cite{Olver:1993} for a proof. On the other hand in order to calculate the corresponding fluxes $g_2$ and $g_4$ we apply the homotopy operator \cite{Olver:1993} to $D_tf_2(\textbf{u})$ and $D_tf_4(\textbf{u})$ respectively. The homotopy operator is given by
$$
H(f(\textbf{u}))=\int_0^1\sum_{i=1}^d\text{I}_i(f)(\lambda \textbf{u})\frac{d\lambda}{\lambda}
$$
where 
$$
\text{I}_{i}(f)=\sum_{k=1}^{\infty}\left(\sum_{s=0}^{k-1}u_{i_{sx}}(-D_x)^{k-s-1}\right)f_{u_{i_{kx}}}.
$$

\end{Proof}

\begin{Cor}
  Let $\rS^1$ be the unitary circle, i.e., $[0,1]$ with matching
  endpoints. Then from Proposition \ref{the:cons} it follows that, upon
  defining
  \begin{equation}
    F_2(\vec u) := \int_{\rS^1} f_2(\vec u) \d x
  \end{equation}
  as the \emph{momentum} functional and
  \begin{equation}
    F_4(\vec u) := \int_{\rS^1} f_4(\vec u) \d x
  \end{equation}
  as the \emph{energy} functional for periodic solutions we have
  \begin{equation}
    \D_t F_2(\vec u) = \D_t F_4(\vec u) = 0.
  \end{equation}
  Moreover this does not just hold for periodic solutions over
  $\rS^1$. Indeed, one can consider the equation (\ref{eq:vmkdv}) over
  $\reals$ and require that solutions decay at infinity, for example
  Schwartz functions, and the result holds. A particular example of
  these are the much celebrated soliton solutions. The conservation
  laws also allow for the a priori control of the solution.
\end{Cor}

\begin{Pro}[Stability bound]
  \label{the:apriori}
  Let the \vmkdv system (\ref{eq:vmkdv}), defined over $\rS^1$, be coupled with initial
  conditions $\vec u_0$ satisfying $F_2(\vec u_0) = C_2 < \infty$ and $F_4(\vec
  u_0) = C_4 < \infty$ then $\vec u$ satisfies
  \begin{equation}
    \Norm{\vec u_x(t)}_{\leb{2}(\rS^1)}
    \leq
    \qp{4C_4
      +
      \frac{C_{GN}^8
      C_2^{3}}{2}}^{1/2}
    ,
  \end{equation}
  where $C_{GN}$ is a constant appearing from the Gagliardo-Nirenberg interpolation inequality.
\end{Pro}
\begin{Proof}
  In view of the definition of $F_4(\vec u)$ we have that
  \begin{equation}
    \label{eq:regpf1}
    \begin{split}
      \Norm{\vec u_x}_{\leb{2}(\rS^1)}^2
      &=
      2 F_4(\vec u)
      +
      \frac 1 4\Norm{\vec u}_{\leb{4}(\rS^1)}^4
      \\
      &=
      2 F_4(\vec u_0)
      +
      \frac 1 4\Norm{\vec u}_{\leb{4}(\rS^1)}^4,
    \end{split}
  \end{equation}
  through the conservativity of $F_4(\vec u)$ given in Theorem
  \ref{the:cons}. Now making use of the Gagliardo-Nirenberg
  interpolation inequality there exists a constant $C_{GN}$ such that
  \begin{equation}
    \Norm{\vec u}_{\leb{4}(\rS^1)}
    \leq
    C_{GN}
    \Norm{\vec u}_{\leb{2}(\rS^1)}^{3/4}
    \Norm{\vec u_x}_{\leb{2}(\rS^1)}^{1/4},
  \end{equation}
  hence
  \begin{equation}
    \label{eq:regpf2}
    \begin{split}
      \frac 1 4\Norm{\vec u}_{\leb{4}(\rS^1)}^4
      &\leq
      \frac 1 4C_{GN}^4
      \Norm{\vec u}_{\leb{2}(\rS^1)}^{3}
      \Norm{\vec u_x}_{\leb{2}(\rS^1)}
      \\
      &\leq 
      \frac{1}{32}
      C_{GN}^8
        {\Norm{\vec u}_{\leb{2}(\rS^1)}^{6}}
        +
      \frac 1 2
        {\Norm{\vec u_x}_{\leb{2}(\rS^1)}^2}
       ,
    \end{split}
  \end{equation}
  through Young's inequality. Substituting (\ref{eq:regpf2}) into (\ref{eq:regpf1}) we see
  \begin{equation}
    \begin{split}
      \frac 1 2 \Norm{\vec u_x}_{\leb{2}(\rS^1)}^2
      &\leq
      2 F_4(\vec u_0)
      +
      \frac{C_{GN}^8}{32}
      \Norm{\vec u}_{\leb{2}(\rS^1)}^{6}
      \\
      &\leq
      2 F_4(\vec u_0)
      +
      \frac{C_{GN}^8}{4}
      F_2(\vec u)^{3}
      \\
      &\leq
      2 F_4(\vec u_0)
      +
      \frac{C_{GN}^8}{4}
      F_2(\vec u_0)^{3}
      \\
      &\leq
      2 C_4
      +
      \frac{C_{GN}^8C_2^{3}}{4}
      ,
    \end{split}
  \end{equation}
  using the conservativity of $F_2(\vec u)$, concluding the proof.
\end{Proof}

\begin{Rem}[Hierarchy of conservation laws]
  Note that the vmKdV equation (\ref{eq:vmkdv})
  admits an infinite hierarchy of conserved quantities. For example,
  after $F_2(\vec{u})$ and $F_4(\vec{u})$ the next member of the
  hierarchy is
  \begin{equation}
    F_6(\vec{u})=\int_{{S}^1}\frac{1}{2}\norm{\vec{u}}^6+10\qp{\vec{u} \cdot \vec{u}_x}^2+\norm{\vec{u}}^2\norm{\vec{u}_x}^2+7\norm{\vec{u}}^2\qp{\vec{u}\cdot \vec{u}_{xx}}+4\norm{\vec{u}_{xx}}^2\d x.
  \end{equation}
  A generating function of the conserved densities for the vmKdV is
  constructed using its Lax representation in
  \cite{AdamapoulouPapamikos:2017}.

  Together with the Gagliardo-Nirenberg interpolation inequality one
  may derive a priori bounds of a similar form to that given in
  Theorem \ref{the:apriori} but in higher order norms. Indeed, for
  $s\in\naturals$ the conservation law $F_{2s}$ naturally gives rise
  to a stability bound in $\sobh{s-1}$.
\end{Rem}

\subsection{Exact solutions to the \vmkdv system}
\label{sec:exact-sol}
The \vmkdv equation (\ref{eq:vmkdv}) is integrable and has already
drawn some attention
\cite{AdamapoulouPapamikos:2017,AncoNgatatWilloughby:2011}.  Its
integrability properties were derived using the structure equation for
the evolution of a curve embedded in an $n$-dimensional Riemannian
manifold with constant curvature
\cite{SandersWang:2003,Anco:2006,Mari-BeffaSandersWang:2002}. The
associated Cauchy problem can be studied analytically using the
inverse scattering transform
\cite{AblowitzClarkson:1991,NovikovManakovPitaevskiiZakharov:1984,FaddeevTakhtajan:2007}. Due
to the fact that it admits a zero curvature representation (or a Lax
representation, see \cite[c.f.]{Lax:1976}), \ie it can be written in
the following form:
\begin{equation}\label{eq:zcc}
  U_t-V_x+\left[U,V\right]=0,
\end{equation}
where $U=U(\vec{u};\lambda)$ and $V=V(\vec{u};\lambda)$ are
appropriate matrices in a Lie algebra having a polynomial dependence
on a spectral parameter $\lambda\in\mathbb{C}$. One can construct, see
\cite{AdamapoulouPapamikos:2017}, a Darboux matrix $M$
\cite{RogersSchief:2002,SalleMatveev:1991} that maps the pair $(U,V)$
to
\begin{equation}\label{eq:DT}
(U,V)\mapsto(\widetilde{U},\widetilde{V})=(MUM^{-1}+M_xM^{-1},MUM^{-1}+M_tM^{-1})
\end{equation}
and $\widetilde{U}=U(\widetilde{\vec{u}};\lambda)$ and
$\widetilde{V}=V(\widetilde{\vec{u}};\lambda)$. In other words
$\widetilde{U}$ and $\widetilde{V}$ have the same structure as $U$ and
$V$ respectively. The transformation \eqref{eq:DT} implies a nonlocal
symmetry $\vec{u} \mapsto \widetilde{\vec{u}}$ of the vmKdV, known as
a B\"acklund transformation. Such transformations that have
applications in geometry \cite{RogersSchief:2002} are characteristic of integrable
equations. Starting with the trivial background solution $\vec{u}=\vec 0$
one can then recursively and algebraically construct the soliton
solutions of vmKdV equation (\ref{eq:vmkdv}). For example, when $d=2$
a 1-soliton solution is given by
\begin{equation}\label{eq:1sol}
\vec{u}=\frac{2\mu}{\cosh\xi_{\mu}}\vec{E},
\end{equation}  
where $\mu\in\mathbb{R}$, $\xi_{\mu}=\mu \qp{x - c_\mu} - \mu^3t$, for
some shift $c_\mu\in\reals$ and $\vec{E}$ is a constant unit vector. A
2-soliton solution is given by
\begin{equation}\label{eq:2sol}
\vec{u}=\frac{F_{\mu,\nu}}{G}\vec{E}_1+\frac{F_{\nu,\mu}}{G}\vec{E}_2,
\end{equation}
where $\vec{E}_1$ and $\vec{E}_2$ are constant unit vectors,
$\mu,~\nu\in\mathbb{R}$ with $\mu\neq\pm \nu$ and
\begin{equation}
  \label{eq:F}
F_{k,l}=2(l^2-k^2)l\cosh\xi_{k}
\end{equation}
and
\begin{equation}
  \label{eq:G}
G=(\mu^2+\nu^2)\cosh\xi_{\mu}\cosh\xi_{\nu}-2\mu\nu \sinh\xi_{\mu}\sinh\xi_{\nu}-2\mu\nu \vec{E}_1\cdot \vec{E}_2.
\end{equation}
The 1-soliton \eqref{eq:1sol} and 2-soliton \eqref{eq:2sol} solutions,
while elegant, are not the most general of their kind, see
\cite{AdamapoulouPapamikos:2017} for details. Nevertheless, the exact
solutions \eqref{eq:1sol} and \eqref{eq:2sol} are both perfectly
adequate for benchmarking our scheme which we shall use them for in
\S\ref{sec:numerics}. Such solutions can also be derived using
Hirota's bilinear form \cite{Hirota:2004,AncoNgatatWilloughby:2011}.

Solitons are, however, a special class of solution for this problem
with a very particular structure. In general one cannot write down
closed form solutions for this problem motivating the need for long
time accurate numerical schemes. We shall proceed by describing the
Hamiltonian structure of the \vmkdv problem which forms the basis for
the design of our numerical scheme.

\begin{Rem}[Problem reformulation and motivation]
  \label{rem:reformulation}
  Since the \vmkdv system is Hamiltonian it can be written as
  \begin{equation}
    \vec u_t = \cP(\vec u) \frac{\delta F_4(\vec u)}{\delta \vec u},
  \end{equation}
  where $\cP(\vec u)$ is a Hamiltonian operator, $F_4(\vec u)$ is the
  induced Hamiltonian and $\frac{\delta \cdot}{\delta \vec u}$ denotes the
  first variation with respect to $\vec u$ \cite{Olver:1993}. For this specific problem the
  Hamiltonian operator acts on a real, $d$-vector function $\vec y$
  and takes the form  \cite{Anco:2006}
  \begin{equation}
    \label{eq:ham-op}
    \begin{split}
      \cP(\vec u) \vec y
      :=
      \vec y_x
      +
      \vec u \bigg \lrcorner \qb{\D_x^{-1} \qp{\vec y \otimes \vec u - \vec u \otimes \vec y}},
    \end{split}
  \end{equation}
  where $\D_x^{-1}$ is the formal inverse operator of $\D_x$,
  $\otimes$ is the tensor product between vectors and
  $\lrcorner$ is an interior product defined through
  \begin{equation}
    \vec x \lrcorner \qp{\vec y \otimes \vec z}
    =
    \qp{\vec x \cdot \vec y}\vec z
    .
  \end{equation}
  This then induces a
  \emph{Poisson bracket}
  \begin{equation}
    \{ F, G \} := \int_{\reals}
    \frac{\delta F}{\delta u} \cP(\vec u) \frac{\delta G}{\delta u} \d x,
  \end{equation}
  a skew-symmetric bilinear form satisfying the Jacobi identity. In
  view of the skew-symmetry of $\cP(\vec u)$ we have
  \begin{equation}
    \D_t F_4(\vec u) = \{ F_4(\vec u), F_4(\vec u) \} = 0.
  \end{equation}
  Notice also that the \vmkdv system can also be written as
  \begin{equation}
    \vec u_t = \{\vec u, F_4(\vec u)\}.
  \end{equation}
  
  The main idea behind the discretisation we propose is to correctly
  \emph{represent the Hamiltonian operator in the finite element space
    whilst preserving the skew-symmetry property of the underlying
    bracket}. Indeed, the proof of Proposition \ref{the:cons} motivates
  rewriting the \vmkdv system by introducing auxiliary variables to
  represent different components of the Hamiltonian operator. We
  consider seeking the tuple $\qp{\vec u,\vec v,\vec w}$ such that
  \begin{equation}
    \label{eq:mixed-system}
    \begin{split}
      \vec 0 &= \vec u_t + \vec v_x + \vec w 
      \\
      \vec 0 &= \vec v - \frac 1 2 \norm{\vec u}^2 \vec u - \vec u_{xx}
      \\
      \vec 0 &= \vec w - \norm{\vec u}^2 \vec u_x + \qp{\vec u_x \cdot \vec u} \vec u.
    \end{split}
  \end{equation}
  Notice that $\vec v = \frac{\delta F_4(\vec u)}{\delta u}$ and $\vec
  w = \vec u \bigg \lrcorner \qb{\D_x^{-1} \qp{\vec v \otimes \vec u -
      \vec u \otimes \vec v}}$. This form of $\vec w$ is extremely
  important as the Hamiltonian operator given in (\ref{eq:ham-op}) is
  nonlocal. The fact that it can be ``localised'' by removing the
  $\D_x^{-1}$ allows for the efficient approximation by Galerkin
  methods.

  This reformulation also means that in the case both arguments of the
  Poisson bracket are the Hamiltonian we may write
  \begin{equation}
    \label{eq:PB}
    0 = \{ F_4(\vec u), F_4(\vec u) \} = \int_\reals \vec v \cdot \qp{\vec v_x + \vec w} \d x.
  \end{equation}
  It is exactly this structure that we try to exploit.
\end{Rem}

\begin{Rem}[Relation to the scalar case]
  As already mentioned when $d=1$, the problem reduces to the mKdV
  equation. In this case $w\equiv 0$ and the mixed system coincides
  with that proposed in \cite{JackamanPryer:2017}. Energy conservative
  schemes can be derived and proven to converge using the same
  techniques presented in that work. For $d>1$, $\vec w$ is
  not necessarily zero and represents the additional contribution arising from
  the Hamiltonian operator described in Remark \ref{rem:reformulation}.
\end{Rem}

\begin{Pro}[The mixed system is conservative]
  \label{pro:mixed-sys-cons}
  Let $\vec u, \vec v, \vec w$ be given by (\ref{eq:mixed-system}) then we have that
  \begin{equation}
    \D_t F_4(\vec u) = \D_t \int_{\rS^1} \frac 1 2 \norm{\vec u_x}^2 - \frac 1 8 \norm{\vec u}^4 \d x = 0.
  \end{equation}
\end{Pro}
\begin{Proof}
  Since the mixed system is equivalent to the \vmkdv system the proof is clear through Proposition \ref{the:cons}, however, for illustrative purposes we present it in full as it will become the basis for the design of our numerical scheme. To begin note
  \begin{equation}
    \begin{split}
      \D_t F_4(\vec u)
      &=
      \int_{\rS^1} \vec u_x \cdot \vec u_{xt} - \frac 1 2 \norm{\vec u}^2 \vec u \cdot \vec u_t \d x
      \\
      &=
      \int_{\rS^1} -\vec u_{xx} \cdot \vec u_{t} - \frac 1 2 \norm{\vec u}^2 \vec u \cdot \vec u_t \d x
      \\
      &=
      \int_{\rS^1} -\vec v \cdot \vec u_{t} \d x.
    \end{split}
  \end{equation}
  Now making use of (\ref{eq:mixed-system})
  \begin{equation}
    \begin{split}
      \D_t F_4(\vec u)
      &=
       \int_{\rS^1} \vec v \cdot \qp{\vec v_x + \vec w} \d x
      \\
      &=
       \int_{\rS^1} \vec v \cdot \vec w \d x
      \\
      &=
       \int_{\rS^1} \qp{\frac 1 2 \norm{\vec u}^2 \vec u + \vec u_{xx}} \cdot \vec w \d x.
    \end{split}
  \end{equation}
  Note that from (\ref{eq:mixed-system}) we can see that $\vec w \cdot \vec u = 0$ and hence
  \begin{equation}
    \label{eq:0}
    \begin{split}
      \D_t F_4(\vec u)
      &=
      \int_{\rS^1} \vec u_{xx} \cdot \vec w \d x
      \\
      &=
      \int_{\rS^1} \vec u_{xx} \cdot \qp{\norm{\vec u}^2 \vec u_x - \qp{\vec u_x \cdot \vec u}\vec u} \d x.
    \end{split}
  \end{equation}
  Now, through an integration by parts we have
  \begin{equation}
    \begin{split}
      \int_{\rS^1} \norm{\vec u}^2 \vec u_x \cdot  \vec u_{xx}
      &=
      - \int_{\rS^1} \qp{\norm{\vec u}^2 \vec u_x}_x \cdot  \vec u_{x}
      \\
      &=
      - \int_{\rS^1} 2 \qp{\vec u \cdot \vec u_x} \qp{\vec u_x \cdot  \vec u_{x}}
      +
      \norm{\vec u}^2 \vec u_{xx} \cdot \vec u_x
      \d x
    \end{split}
  \end{equation}
  and hence
  \begin{equation}
    \label{eq:1}
    \int_{\rS^1} \norm{\vec u}^2 \vec u_x \cdot  \vec u_{xx}
    =
    -\int_{\rS^1} \qp{\vec u \cdot \vec u_x} \qp{\vec u_x \cdot  \vec u_{x}}.
  \end{equation}
  In addition,
  \begin{equation}
    \begin{split}
      \int_{\rS^1} \qp{\vec u_x \cdot \vec u} \qp{\vec u\cdot \vec u_{xx}} \d x
      &=
      -\int_{\rS^1} \qp{{\vec u_x \cdot \vec u} \vec u}_x \cdot \vec u_{x} \d x
      \\
      &=
      -\int_{\rS^1} \qp{\vec u_{xx} \cdot \vec u} \qp{\vec u \cdot \vec u_{x}}
      +
      2\qp{\vec u_{x} \cdot \vec u_x} \qp{\vec u \cdot \vec u_{x}}
      \d x
    \end{split}
  \end{equation}
  and hence
  \begin{equation}
    \label{eq:2}
    \int_{\rS^1} \qp{\vec u_x \cdot \vec u} \qp{\vec u\cdot \vec u_{xx}} \d x
    =
    -\int_{\rS^1} \qp{\vec u_{x} \cdot \vec u_x} \qp{\vec u \cdot \vec u_{x}}
      \d x.
  \end{equation}
  Substituting (\ref{eq:1}) and (\ref{eq:2}) into (\ref{eq:0}) concludes the proof.
\end{Proof}

\section{Temporal discretisation}
\label{sec:temp}

For the readers convenience we will present an argument for designing
the temporally discrete {scheme} in the spatially continuous setting.
We consider a time interval $[0,T]$ subdivided into a partition of $N$
consecutive adjacent subintervals whose endpoints are denoted
$t_0=0<t_1<\ldots<t_{N}=T$.  The $n$-th timestep is defined as
${\tau_n := t_{n+1} - t_{n}}$.  We will consistently use the shorthand
$y^n(\cdot):=y(\cdot,t_n)$ for a generic time function $y$. We also
denote $\nplush{y} := \frac{1}{2}\qp{y^n + y^{n+1}}$.

We consider the temporal discretisation of the mixed system (\ref{eq:mixed-system}) as follows: Given $\vec u^0$, for $n\in [0, N]$ find $\vec u^{n+1}$ such that
\begin{equation}
  \label{eq:temp-mixed-system}
  \begin{split}
    \vec 0 &= \frac{\vec u^{n+1} - \vec u^n}{\tau_n} + \vec v_x^{n+1} + \vec w^{n+1} 
    \\
    \vec 0 &= \vec v^{n+1} - \frac 1 2 \qp{\norm{\vec u^n}^2 + \norm{\vec u^{n+1}}^2} \nplush{\vec u} - \nplush{\vec u_{xx}}
    \\
    \vec 0 &= \vec w^{n+1} - \norm{\nplush{\vec u}}^2 \nplush{\vec u_x} + \qp{\nplush{\vec u_x} \cdot \nplush{\vec u}} \nplush{\vec u}.
  \end{split}
\end{equation}

\begin{Rem}[Structure of the temporal discretisation]
  The temporal discretisation given in (\ref{eq:temp-mixed-system}) is
  \emph{not} a Runge-Kutta method. It resembles a Crank-Nicolson
  discretisation, however the treatment of the nonlinearity is
  different. It is formally of second order and is constructed such
  that it satisfies the next Theorem. Although construction of
  higher order methods is possible they become very complicated to
  write down so we will not press this point here.
\end{Rem}

\begin{The}[Conservativity of the temporal discretisation]
  Let $\{ \vec u^n \}_{n=0}^N$ be a temporally discrete solution of
  (\ref{eq:temp-mixed-system}) then we have
  \begin{equation}
    F_4(\vec u^{n}) = F_4(\vec u^0) \Foreach n \in [0, N].
  \end{equation}
\end{The}
\begin{Proof}
  It suffices to show that 
  \begin{equation}
    F_4(\vec u^{n+1}) - F_4(\vec u^n) = 0
  \end{equation}
  and then the result follows inductively. So
  \begin{equation}
    \begin{split}
      2\qp{F_4(\vec u^{n+1}) - F_4(\vec u^n)}
      &=
      \int_{\rS^1} \norm{\vec u_x^{n+1}}^2
      -
      \norm{\vec u_x^{n}}^2
      -
      \frac 1 4 \norm{\vec u^{n+1}}^4
      +
      \frac 1 4 \norm{\vec u^{n}}^4 \d x
      \\
      &=
      \int_{\rS^1}
      \qp{\vec u_x^{n+1} - \vec u_x^{n}} \cdot \qp{\vec u_x^{n+1} + \vec u_x^{n}}
      \\
      &\qquad -
      \frac 1 4 \qp{\vec u^{n+1} - \vec u^{n}} \cdot
      \qp{
        \norm{\vec u^{n+1}}^2 \vec u^{n+1}
        +
        \norm{\vec u^{n+1}}^2 \vec u^{n}
        +
        \norm{\vec u^{n}}^2 \vec u^{n+1}
        +
        \norm{\vec u^{n}}^2 \vec u^{n}
      }
      \d x
      \\
      &=
      -\int_{\rS^1}
      \qp{\vec u^{n+1} - \vec u^{n}} \cdot \qp{\vec u_{xx}^{n+1} + \vec u_{xx}^{n}}
      \\
      &\qquad -
      \frac 1 2 \qp{\vec u^{n+1} - \vec u^{n}} \cdot
      \qp{
        \norm{\vec u^{n+1}}^2 \nplush{\vec u}
        +
        \norm{\vec u^{n}}^2 \nplush{\vec u}
      }
      \d x
      \\
      &=
      -\int_{\rS^1}
      \qp{\vec u^{n+1} - \vec u^{n}} \cdot \vec v^{n+1} \d x,
    \end{split}
  \end{equation}
  through expanding differences, integrating by parts and using the
  scheme (\ref{eq:temp-mixed-system}). Now, again using the scheme
  \begin{equation}
    \begin{split}
      2\qp{F_4(\vec u^{n+1}) - F_4(\vec u^n)}
      &=
      \int_{\rS^1}
      \tau_n \qp{\vec v_x^{n+1} + \vec w^{n+1}} \cdot \vec v^{n+1} \d x
      \\
      &=
      \int_{\rS^1}
      \tau_n \vec w^{n+1} \cdot \vec v^{n+1} \d x
      \\
      &=
      \int_{\rS^1}
      \tau_n \vec w^{n+1} \cdot
      \qp{
        \frac 1 2
        \qp{\norm{\vec u^{n+1}}^2 \nplush{\vec u}
          +
          \norm{\vec u^{n}}^2 \nplush{\vec u}
        }
        +
        \nplush{\vec u_{xx}}
      }
      \d x
      \\
      &=
      \int_{\rS^1}
      \tau_n \vec w^{n+1} \cdot
      \nplush{\vec u_{xx}}
      \d x,
    \end{split}
  \end{equation}
  in view of the orthogonality condition $\vec w^{n+1}\cdot
  \nplush{\vec u} = 0$ following from the third equation of
  (\ref{eq:temp-mixed-system}). Now we may use the definition of
  $\vec w^{n+1}$ and the identities (\ref{eq:1}) and (\ref{eq:2})
  to conclude.
\end{Proof}

\begin{Rem}[Conservation of other invariants]
  \label{rem:cons}
  This discretisation does not lend itself to conservation of other
  invariants, for example even the quadratic invariant $F_2$ is not
  conserved under this scheme. A class of Runge-Kutta methods that are
  able to exactly conserve all quadratic invariants are the
  Gauss-Radau family, this is because they are \emph{symplectic}. When
  one considers higher order invariants, it seems that schemes must be
  designed individually and there is no class that can exactly
  conserve all. 
\end{Rem}

\section{Spatial and full discretisation}
\label{sec:spat}

In this section we describe the discretisation which we analyse for
the approximation of (\ref{eq:vmkdv}). We show that the scheme has a
constant energy functional consistent with that of the original PDE
system.
\begin{Defn}[Finite element space]
  We discretise (\ref{eq:vmkdv}) spatially using a piecewise polynomial
  continuous finite element method. To that end we let $\rS^1 := [0,1]$
  be the unit interval with matching endpoints and choose
  \begin{equation}
    0 =
    x_0 < x_1 < \dots < x_M
    = 1.
  \end{equation}
  Note that in the numerical experiments we take a larger periodic
  interval, however for clarity of presentation we restrict our
  attention in this section to $\rS^1$. We denote $K_m=[x_m,x_{m+1}]$
  to be the $m$--th subinterval and let $h_m:= \norm{K_m}$ be its
  length with $\T{} = \{ K_m \}_{m=0}^{M-1}$. We impose that the ratio
  $h_m/h_{m+1}$ is bounded from above and below for
  $m=0,\dots,M-1$. Let $\poly q$ be the space of polynomials of degree
  less than or equal to $q$, then we introduce the \emph{finite
    element space}
  \begin{equation}
    \fes 
    :=
    \ensemble{\Phi}
             {\Phi \vert _{K_m} \in \poly q(K_m)\cap \cont{0}(\rS^1)}.
  \end{equation}
\end{Defn}
Throughout this section and in the sequel, we will use capital Latin
letters to denote spatially discrete trial functions and capital Greek
letters to denote discrete test functions.

\subsection{Spatial discretisation}

Before we give the discretisation we propose let us first consider a
direct semi-discretisation of the mixed system (\ref{eq:mixed-system}), to
find $\qp{\vec U, \vec V, \vec W} \in \fes^d \times \fes^d
\times \fes^d$ such that
\begin{equation}
  \label{eq:dis-mixed-system}
  \begin{split}
    0 &= \int_{\rS^1} \qp{\vec U_{t} + \vec V_{x} + \vec W} \cdot \vec \Phi \d x
    \\
    0 &= \int_{\rS^1} \qp{\vec V - \frac 1 2 \norm{\vec U}^2 \vec U} \cdot \vec \Psi
    + \vec U_{x} \cdot \vec \Psi_x \d x
    \\
    0 &= \int_{\rS^1} \qp{\vec W - \norm{\vec U}^2 \vec U_x + \qp{\vec U_x \cdot \vec U} \vec U} \cdot \vec \Xi \d x \Foreach \qp{\vec \Phi, \vec \Psi, \vec \Xi} \in \fes^d \times \fes^d \times \fes^d.
  \end{split}
\end{equation}
One may run through the calculation in the Proof of Proposition
\ref{pro:mixed-sys-cons} analogously to see that
\begin{equation}
  \label{eq:vw}
  \D_t F_4(\vec U) = \int_{\rS^1} \vec V \cdot \vec W \d x,
\end{equation}
whereby in the continuous case one uses the fact that $\vec v$ and
$\vec w$ are orthogonal. In the discrete setting there is no reason
why this should be the case and, indeed, except in very special cases,
it is not. This necessitates a formulation that forces $\int_{\rS^1}
\vec V \cdot \vec W \d x = 0$ thus ensuring conservation of $F_4(\vec
U)$. We achieve this through a Lagrange multiplier approach
encapsulated by the following spatially discrete scheme, to seek
$\qp{\vec U, \vec V, \vec W, P} \in \fes^d \times \fes^d \times \fes^d
\times \reals \slash \{ 0\}$ such that
\begin{equation}
  \label{eq:dis-mixed-system-sd}
  \begin{split}
    0 &= \int_{\rS^1} \qp{{\vec U_t}  + \vec V_{x} + \vec W } \cdot \vec \Phi \d x
    \\
    0 &= \int_{\rS^1} \qp{\vec V - \frac 1 2 \norm{\vec U}^2 \vec U } \cdot \vec \Psi
    + \vec U_{x} \cdot \vec \Psi_x \d x
    \\
    0 &= \int_{\rS^1} \qp{\vec W - \norm{\vec U}^2 \vec U_x + \qp{\vec U_x \cdot \vec U} \vec U} \cdot \vec \Xi \d x 
    \\
    0 &= \int_{\rS^1} P \vec V \cdot \vec \Xi + \vec V \cdot \vec W \Pi  \d x
    \\
    \vec U^0 &= I^0 u_0, \Foreach \qp{\vec \Phi, \vec \Psi, \vec \Xi, \Pi} \in {\fes}^d \times {\fes}^d \times {\fes}^d \times \reals \slash \{ 0\}
  \end{split}
\end{equation}

\begin{The}[Conservativity of the spatially discrete scheme]
  \label{the:cons-sd}
  Let $(\vec U,\vec V,\vec W,P)$ solve the spatially discrete formulation (\ref{eq:dis-mixed-system-sd}) then
  \begin{equation}
    \D_t F_4(\vec U) = \D_t \int_{\rS^1} \frac 1 2 \norm{\vec U_x}^2 - \frac 1 8 \norm{\vec U}^4 \d x = 0.
  \end{equation}
\end{The}

\begin{Proof}
  An analogous argument to the Proof of Proposition
  \ref{pro:mixed-sys-cons} yields (\ref{eq:vw}). To conclude pick $\Pi
  = P$ and $\vec \Xi = \vec W$ to see that
  \begin{equation}
    2P \int_{\rS^1} \vec V \cdot \vec W \d x = 0,
  \end{equation}
  as required.
\end{Proof}

\begin{Rem}[Compatibility of the scheme with the Poisson bracket]
  Notice that in view of Theorem \ref{the:cons-sd} the spatial
  discretisation is compatible with the Poisson structure of the
  \vmkdv system, indeed, using the same formulation as (\ref{eq:PB})
  we have that the numerical scheme can be written as
  \begin{equation}
    \vec U_t = \{ \vec U, \vec F_4(\vec U) \} = -\qp{\pi\qp{\vec V_x} + \vec W},
  \end{equation}
  where $\pi$ denotes the $\leb{2}$ orthogonal projector onto $\fes$.
  In addition, the evolution of the Hamiltonian can be described
  consistently
  \begin{equation}
    \D_t F_4(\vec U) = \{ F_4(\vec U), F_4(\vec U) \} = 0.
  \end{equation}
  It is important to note that the evolution of other quantities, such
  as further invariants are not compatible with this structure, for example
  \begin{equation}
    \D_t F_2(\vec U) = \{ F_2(\vec U), F_4(\vec U) \} \neq 0.
  \end{equation}
\end{Rem}

\subsection{Fully discrete scheme}
Making use of the semi discretisations developed in \S\ref{sec:temp}
and \S\ref{sec:spat} we consider a fully discrete approximation that
consists of finding a sequence of functions
\newline $\qp{\vec U^{n+1}, \vec V^{n+1}, \vec W^{n+1}, P^{n+1}}
\in \qb{\fes^n}^d\times \qb{\fes^n}^d\times \qb{\fes^n}^d \times \reals\slash \{0\}$ such that
for each $n \in [0, N-1]$ we have
\begin{equation}
  \label{eq:dis-mixed-system-fd}
  \begin{split}
    0 &= \int_{\rS^1} \qp{\frac{\vec U^{n+1} - \vec U^n}{\tau_n}  + \vec V_{x}^{n+1} + \vec W^{n+1} } \cdot \vec \Phi \d x
    \\
    0 &= \int_{\rS^1} \qp{\vec V^{n+1} - \frac 1 2 \qp{\norm{\vec U^n}^2 + \norm{\vec U^{n+1}}^2} \vec U^{n+1/2} } \cdot \vec \Psi
    + \vec U_{x}^{n+1/2} \cdot \vec \Psi_x \d x
    \\
    0 &= \int_{\rS^1} \qp{\vec W^{n+1} - \norm{\vec U^{n+1/2}}^2 \vec U_x^{n+1/2} + \qp{\vec U_x^{n+1/2} \cdot \vec U^{n+1/2}} \vec U^{n+1/2}} \cdot \vec \Xi \d x 
    \\
    0 &= \int_{\rS^1} P^{n+1} \vec V^{n+1} \cdot \vec \Xi + \vec V^{n+1} \cdot \vec W^{n+1} \Pi  \d x
    \\
    \vec U^0 &= \pi^0 u_0, \Foreach \qp{\vec \Phi, \vec \Psi, \vec \Xi, \Pi} \in \qb{\fes^n}^d \times \qb{\fes^n}^d \times \qb{\fes^n}^d \times \reals\slash \{0\}
  \end{split}
\end{equation}
where $\pi^0$ denotes the $\leb{2}$ orthogonal projector into
$\fes^0$. This is the direct discretisation of the mixed system
(\ref{eq:mixed-system}) with the temporal discretisation as that
proposed in \S\ref{sec:temp} with an additional equation for a unknown
real number that represents a Lagrange multiplier ensuring
$\int_{\rS^1} \vec V^n \cdot \vec W^n \d x = 0$ for all $n$.

\begin{The}[Conservativity of the fully discrete scheme]
  \label{the:cons-fullydiscrete}
  Let $\{ \vec U^n \}_{n=0}^N$ be the fully discrete scheme generated
  by (\ref{eq:dis-mixed-system-fd}), then we have that
  \begin{equation}
    F_4(\vec U^{n}) = F_4(\vec U^0) \Foreach n \in [0, N].
  \end{equation}
\end{The}
\begin{Proof}
  It suffices to show that 
  \begin{equation}
    F_4(\vec U^{n+1}) - F_4(\vec U^n) = 0
  \end{equation}
  and then the result follows inductively. To this end
  \begin{equation}
    \begin{split}
      2\qp{F_4(\vec U^{n+1}) - F_4(\vec U^n)}
      &=
      \int_{\rS^1} \norm{\vec U_x^{n+1}}^2
      -
      \norm{\vec U_x^{n}}^2
      -
      \frac 1 4 \norm{\vec U^{n+1}}^4
      +
      \frac 1 4 \norm{\vec U^{n}}^4 \d x
      \\
      &=
      \int_{\rS^1}
      \qp{\vec U_x^{n+1} - \vec U_x^{n}} \cdot \qp{\vec U_x^{n+1} + \vec U_x^{n}}
      \\
      &\qquad -
      \frac 1 4 \qp{\vec U^{n+1} - \vec U^{n}} \cdot
      \qp{
        \norm{\vec U^{n+1}}^2 \vec U^{n+1}
        +
        \norm{\vec U^{n+1}}^2 \vec U^{n}
        +
        \norm{\vec U^{n}}^2 \vec U^{n+1}
        +
        \norm{\vec U^{n}}^2 \vec U^{n}
      }
      \d x
      \\
      &=
      -\int_{\rS^1}
      \qp{\vec U^{n+1} - \vec U^{n}} \cdot {\vec V^{n+1}} \d x,
    \end{split}
  \end{equation}
  through expanding differences and using the second equation of (\ref{eq:dis-mixed-system-fd}). Now, using the first equation of (\ref{eq:dis-mixed-system-fd})
  \begin{equation}
    \begin{split}
      2\qp{F_4(\vec U^{n+1}) - F_4(\vec U^n)}
      &=
      \int_{\rS^1}
      \tau_n \qp{\vec V_x^{n+1} + \vec W^{n+1}} \cdot \vec V^{n+1} \d x
      \\
      &=
      \int_{\rS^1}
      \tau_n \vec W^{n+1} \cdot \vec V^{n+1} \d x
      \\
      &=
      0
    \end{split}
  \end{equation}
  using the fourth equation of (\ref{eq:dis-mixed-system-fd}) with
  $\Pi=P^{n+1}$ and $\vec \Xi = \vec W^{n+1}$, concluding the proof.
\end{Proof}

\section{Numerical experiments}
\label{sec:numerics}

In this section we illustrate the performance of the method proposed
through a series of numerical experiments. The brunt of the
computational work was carried out using Firedrake
\cite{Firedrake:2015} and Paraview was used as a visualisation
tool. The code written for this purpose is freely available at
\cite{vmkdvzenodo}. We ignore the effect of numerical integration in
all our computations by taking a sufficiently high quadrature degree
that allows for accurate evaluation of all integrals in all our
numerical examples. For each benchmark test we fix the polynomial
degree $q$ and compute a sequence of solutions with $h = h(i) =
2^{-i}$ and $\tau$ chosen either so $ \tau \ll h$, to make the
temporal discretisation error negligible, or so $\tau = h$ so temporal
discretisation error dominates. This is done for a sequence of
refinement levels, $i=l, \dots, L$.

\begin{Defn}[Experimental order of convergence]
  Given two sequences $a(i)$  and $h(i) \downto 0$ we define the \emph{experimental order of convergence} (EOC) to be the local slope of the $\log{a(i)}$ vs. $\log{h(i)}$ curve, i.e., 
  \begin{equation}
    \EOC(a,h;i) = \frac{\log\qp{\frac{a_{i+1}}{a_i}}}{\log\qp{\frac{h_{i+1}}{h_i}}}.
  \end{equation}
\end{Defn}

\subsection{Test 1 - Asymptotic benchmarking of a $1$-soliton solution}

\begin{figure}[h!]
 \centering
 \caption{
   \label{fig:1soldev}
   Here we examine the conservative discretisation scheme with various polynomial degrees, $q$, approximating the exact solution (\ref{eq:1sol}) with initial conditions given by (\ref{eq:1solIC}). We show the deviation in the two invariants $F_i$, $i=2,4$, corresponding to momentum and energy. In each test we take a fixed spatial discretisation parameter of $h=0.25$ and fixed time step of $\tau = 0.001$. Notice that in each case the deviation in energy is smaller than the solver tolerance of $10^{-12}$ and the deviation in momentum is bounded. In addition, as the degree of approximation is increased the deviation in momentum becomes smaller, in this case by around two orders of magnitude per polynomial order. The simulations are simulated for long time to test conservativity with $T=100$ in each case.}
   \subfigure[][Here $q=1$.]{
   \includegraphics[scale=\figscale, width=0.31\figwidth]{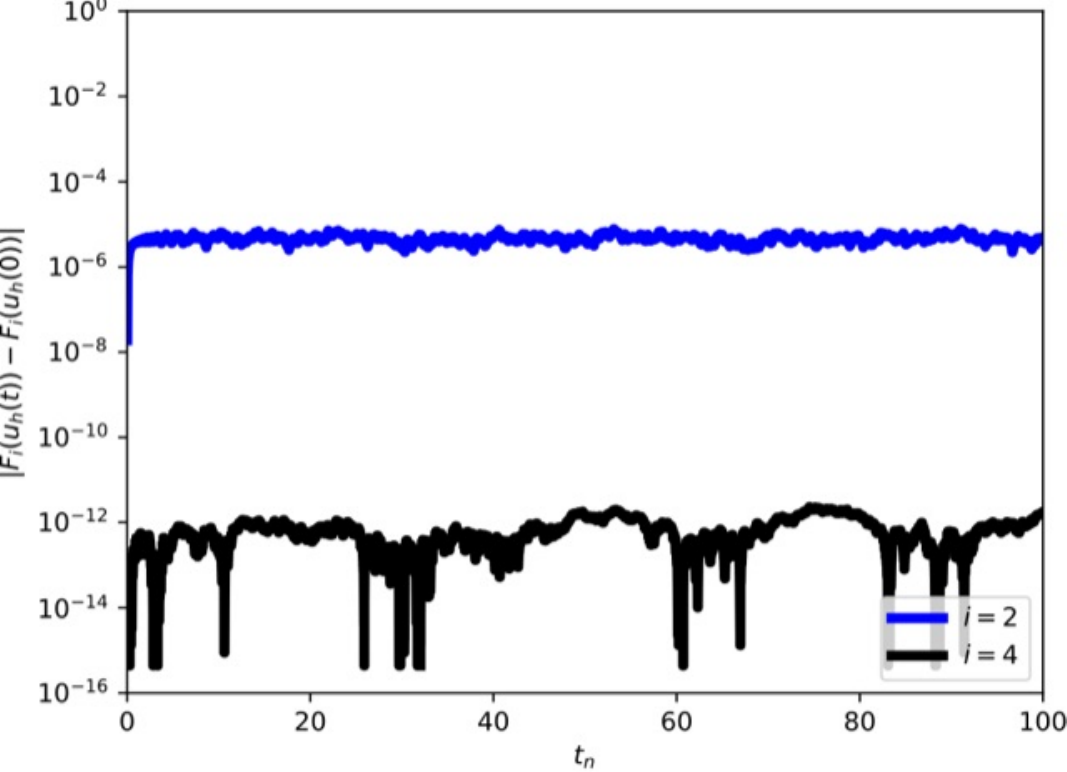}
 }
   \subfigure[][Here $q=2$.]{
   \includegraphics[scale=\figscale, width=0.31\figwidth]{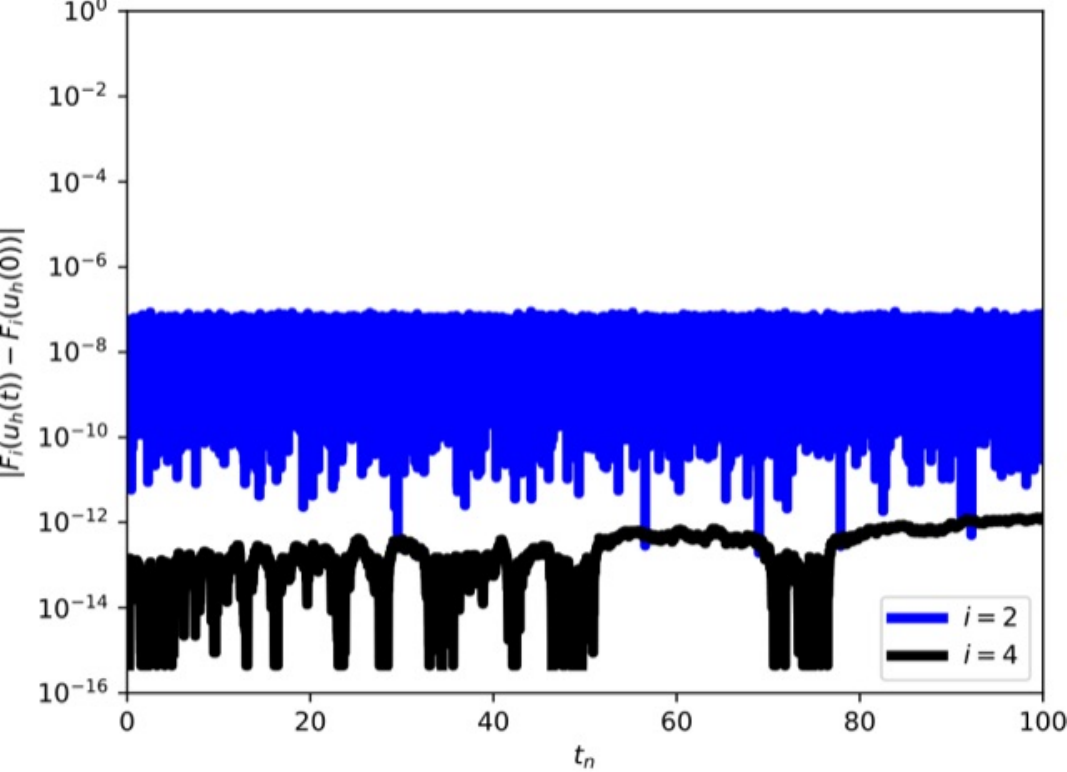}
 } 
  \subfigure[][Here $q=3$.]{
   \includegraphics[scale=\figscale, width=0.31\figwidth]{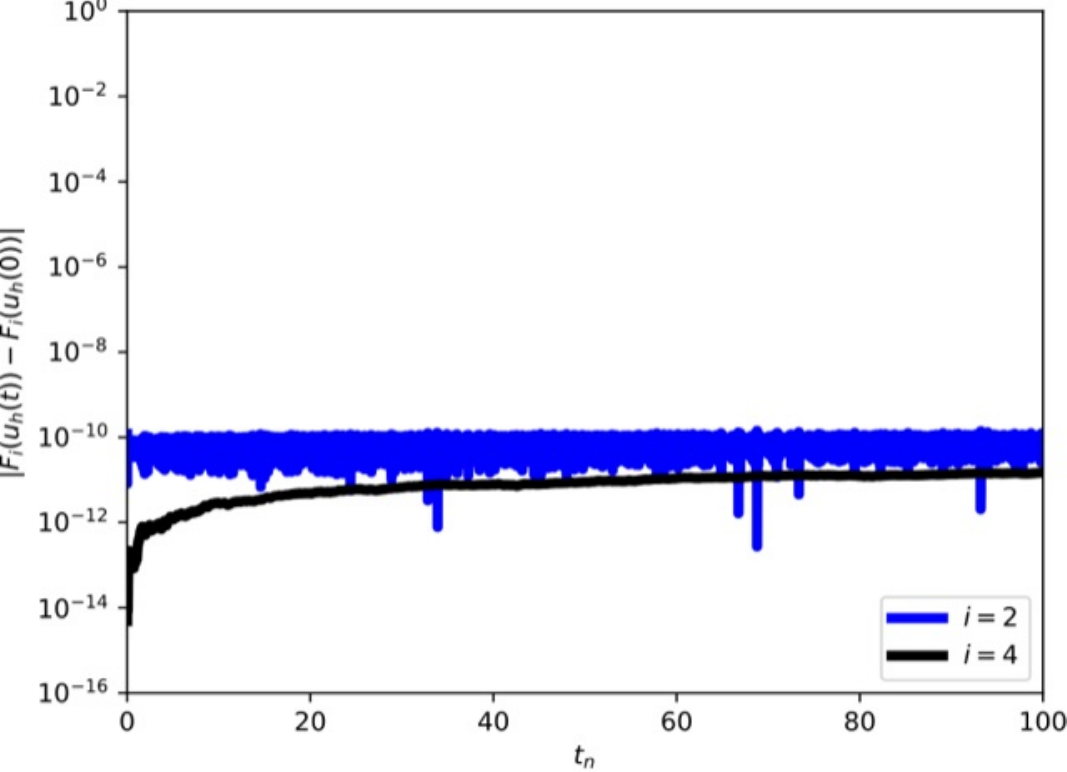}
 }
\end{figure}

\begin{figure}[h!]
 \centering
 \caption{
   \label{fig:1solcon}
   Here we examine the conservative discretisation scheme with various
   polynomial degrees, $q$, approximating the exact solution
   (\ref{eq:1sol}) with initial conditions given by
   (\ref{eq:1solIC}). We show the errors measured in the
   $\leb{\infty}(0,t_n; \leb{2}(\rS^1))$ norm for each component of
   the system and the EOC for test runs that benchmark both the
   spatial and temporal discretisation and show that the scheme is of
   optimal order. We use $e_{u_i}:= \Norm{u_i -
     U_{i}}_{\leb{\infty}(0,t_n; \leb{2}(\rS^1))}$ for $i=1,2$, the
   components of the solution $\vec u = \qp{u_1,u_2}^\transpose$ and
   numerical approximation $\vec U = \qp{U_1,U_2}^\transpose$.}
 
   \subfigure[][Here $q=1$ and we fix $\tau{} = 0.00001$. This is sufficiently small that the spatial discretisation error dominates.]{
   \includegraphics[scale=\figscale, width=0.47\figwidth]{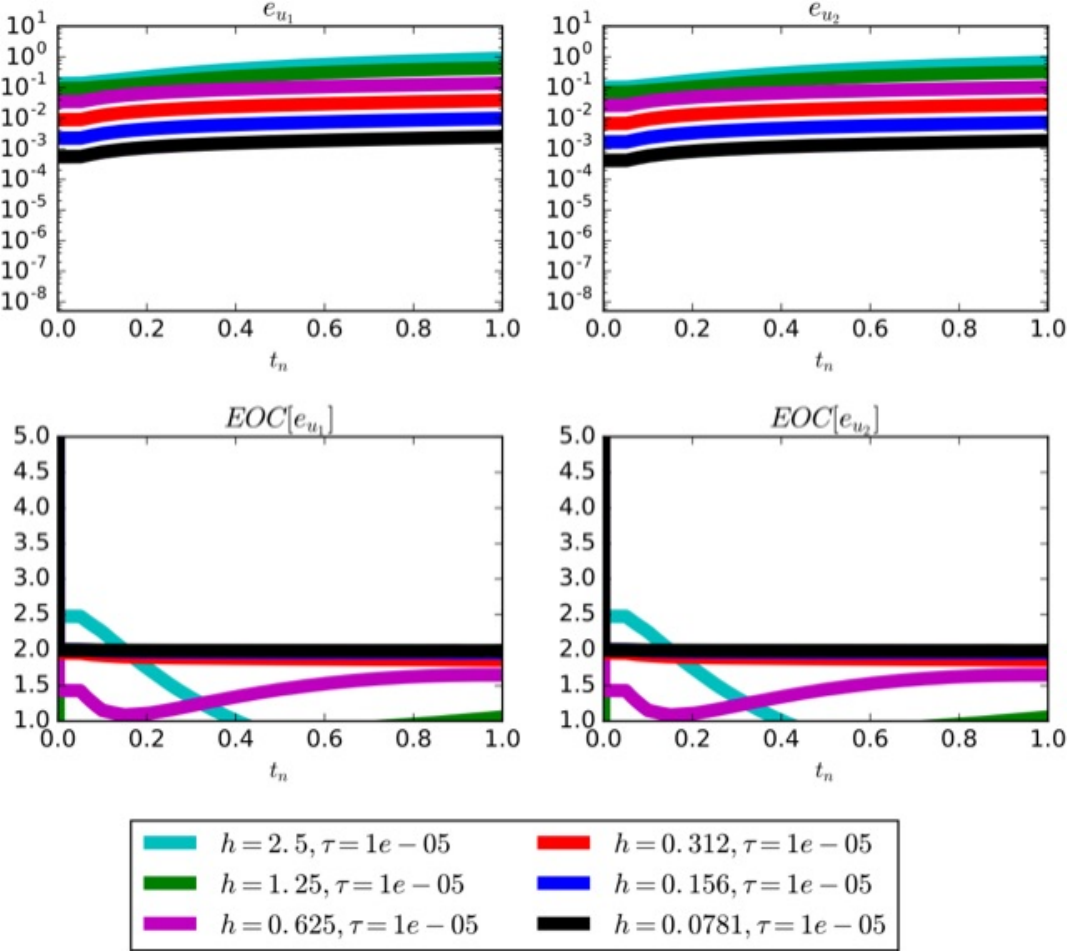}
 }
   \subfigure[][Here $q=2$ and we fix $\tau = 0.00001$. This is sufficiently small that the spatial discretisation error dominates.]{
   \includegraphics[scale=\figscale, width=0.47\figwidth]{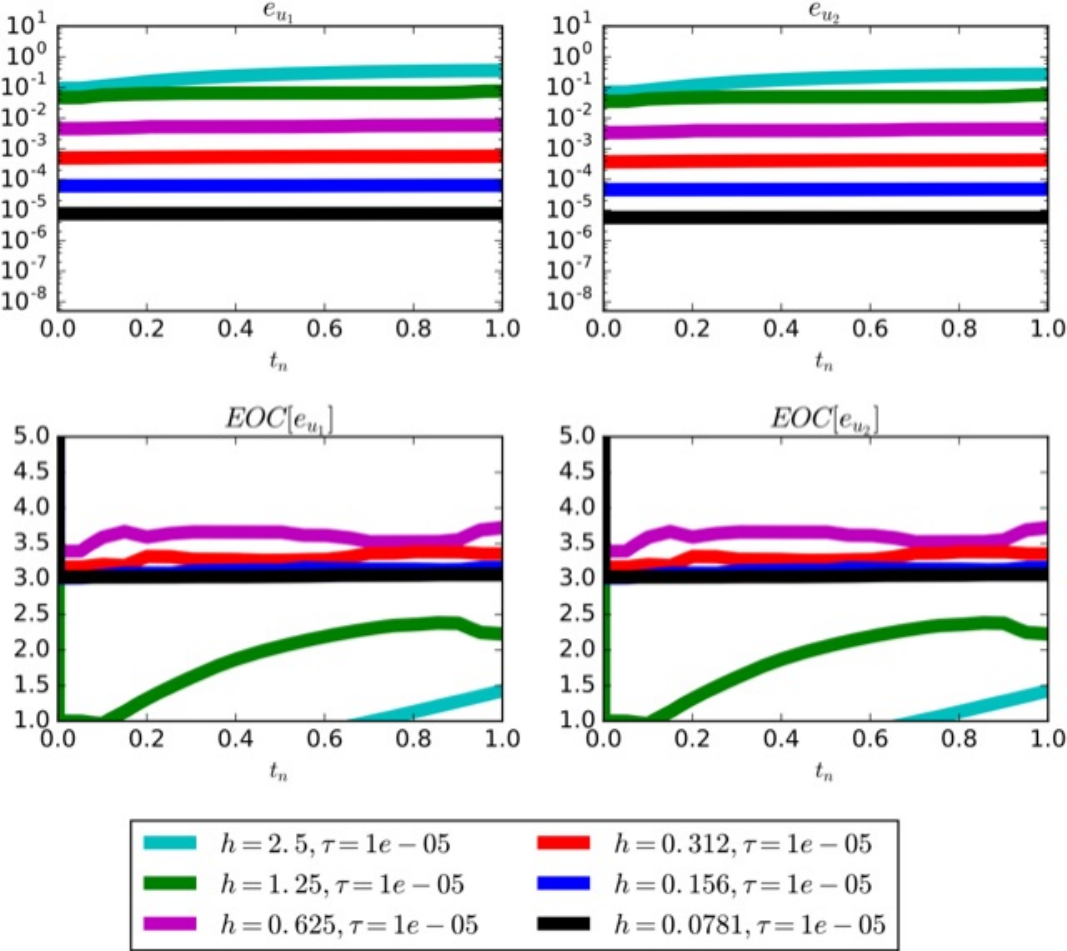}
 } 
  \subfigure[][Here $q=3$ and we fix $\tau = 0.00001$. This is sufficiently small that the spatial discretisation error dominates.]{
   \includegraphics[scale=\figscale, width=0.47\figwidth]{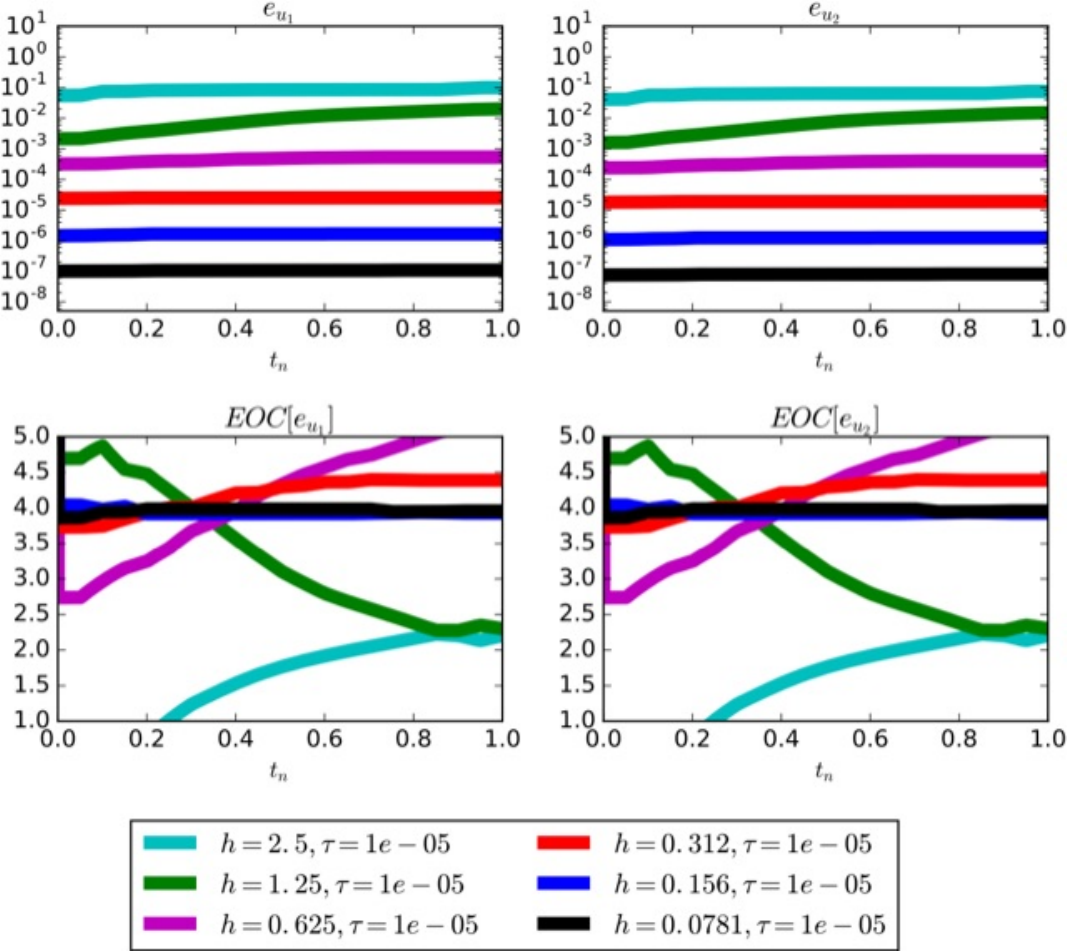}
  }
  \subfigure[][Here $q=2$ and on every refinement level we choose a coupling $\tau = C h$. Note that the time discretisation error here dominates.]{
    \includegraphics[scale=\figscale, width=0.47\figwidth]{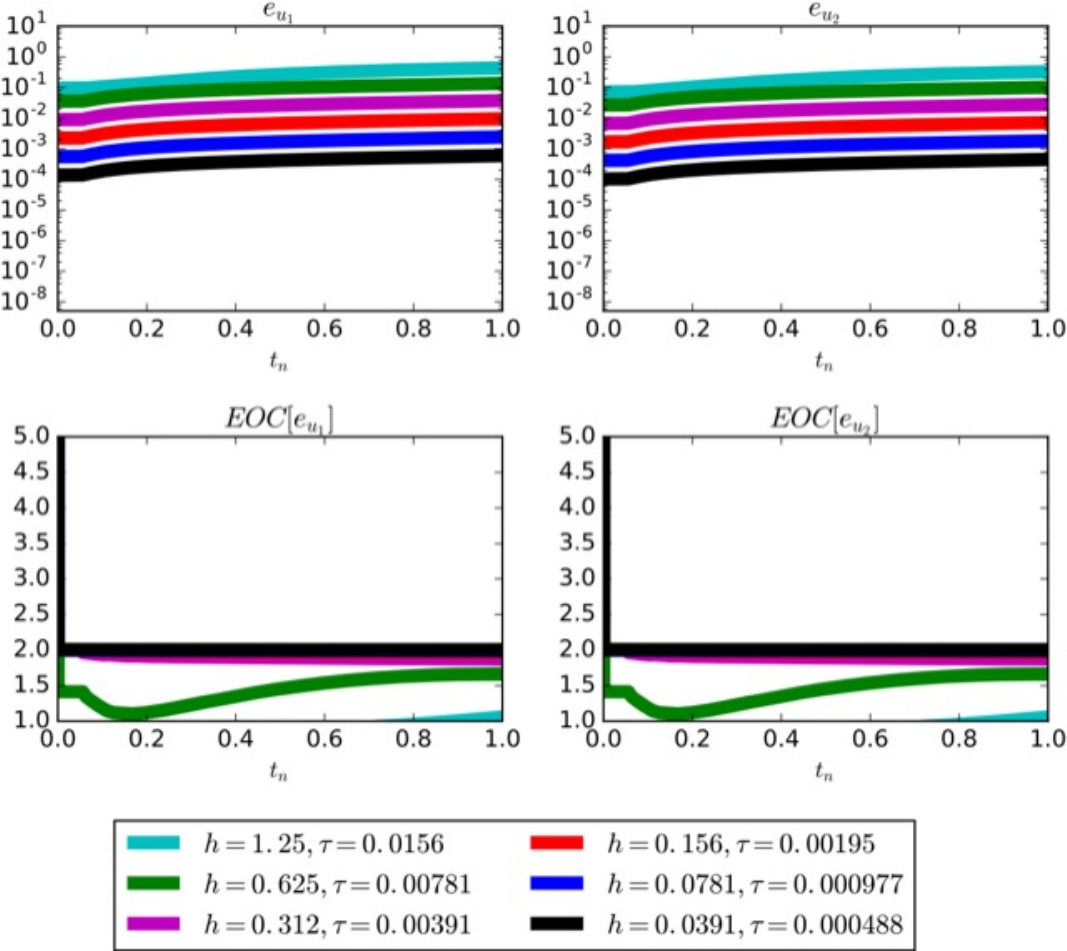}
  }
\end{figure}

We take $d=2$ and
\begin{equation}
  \label{eq:1solIC}
  \vec{u_0}=\frac{2\mu}{\cosh\qp{\mu \qp{x - c_\mu}}}\vec{E},
\end{equation}
over the periodic domain $[0,40]$ with $\vec E = (0.8,
(1-0.8^2)^{0.5})^\transpose$, $\mu=1$ and $c_\mu=20$. The
exact solution is then given by (\ref{eq:1sol}). We take a uniform
timestep and uniform meshes that are fixed with respect to
time. Convergence results are shown in Figure \ref{fig:1solcon} and
conservativity over long time is given in Figure
\ref{fig:1soldev}. Note that for the $1$-soliton solution we have
$\vec W^n \equiv \vec 0$ for all $n$ in which case the Lagrange
multiplier is not required as $\int_{\rS^1} \vec V^n \cdot \vec W^n =
0$ trivially for all $n$.

\clearpage

\subsection{Test 2 - Asymptotic benchmarking of a $2$-soliton solution}

\begin{figure}[h]
 \centering
 \caption{
   \label{fig:2soldev}
      Here we examine the conservative discretisation scheme with various polynomial degrees, $q$, approximating the exact solution (\ref{eq:2sol}) with initial conditions given by (\ref{eq:2solIC}). We show the deviation in the two invariants $F_i$, $i=2,4$, corresponding to momentum and energy respectively. In each test we take a fixed spatial discretisation parameter of $h=0.25$ and fixed time step of $\tau = 0.001$. Notice that in each case the deviation in energy is smaller than the solver tolerance of $10^{-12}$ and the deviation in momentum is bounded. The simulations are simulated for long time to test conservativity with $T=100$ in each case.}
   \subfigure[][Here $q=1$.]{
   \includegraphics[scale=\figscale, width=0.31\figwidth]{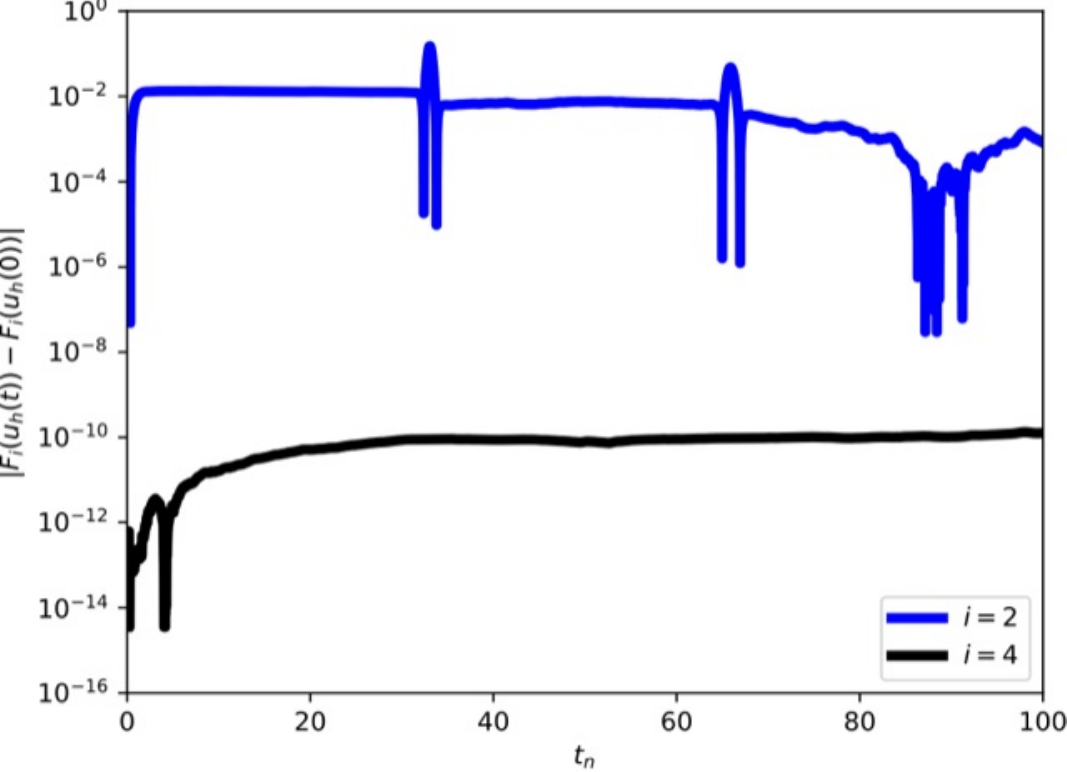}
 }
   \subfigure[][Here $q=2$.]{
   \includegraphics[scale=\figscale, width=0.31\figwidth]{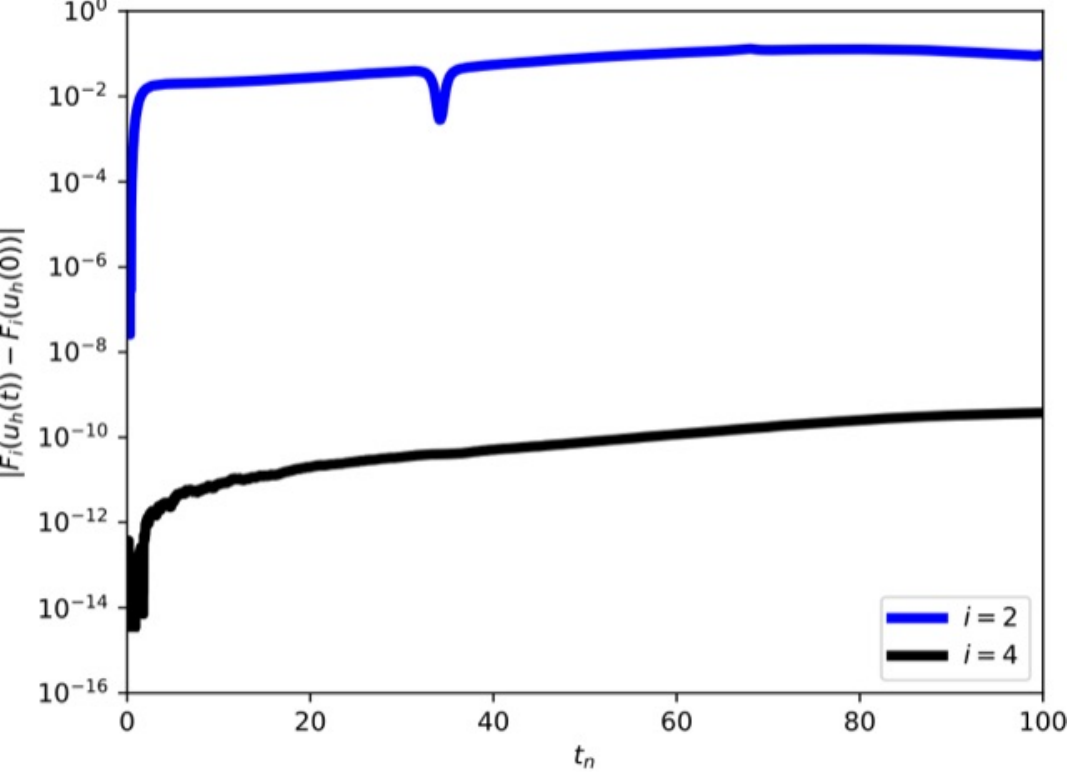}
 } 
  \subfigure[][Here $q=3$.]{
   \includegraphics[scale=\figscale, width=0.31\figwidth]{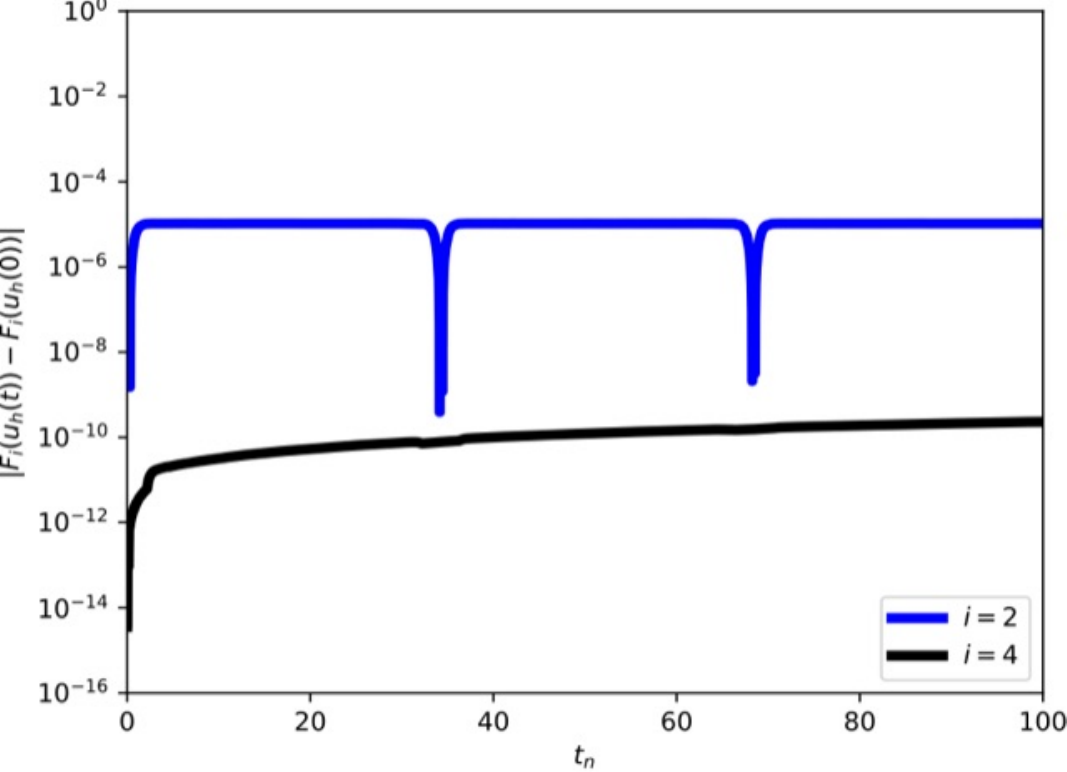}
 }
\end{figure}

\begin{figure}[h]
 \centering
 \caption{
   \label{fig:2solcon}
      Here we examine the conservative discretisation scheme with various polynomial degrees, $q$, approximating the exact solution (\ref{eq:2sol}) with initial conditions given by (\ref{eq:2solIC}). We show the errors measured in the $\leb{\infty}(0,t_n; \leb{2}(\rS^1))$ norm for each component of the system and the EOC for test runs that benchmark both the spatial and temporal discretisation and show that the scheme is of the correct order. We use $e_{u_i}:= \Norm{u_i - U_i}_{\leb{\infty}(0,t_n; \leb{2}(\rS^1))}$ for $i=1,2$, the components of the solution $\vec u = \qp{u_1,u_2}^\transpose$.}
   \subfigure[][Here $q=1$ and we fix $\tau{} = 0.00001$. This is sufficiently small that the spatial discretisation error dominates.]{
   \includegraphics[scale=\figscale, width=0.47\figwidth]{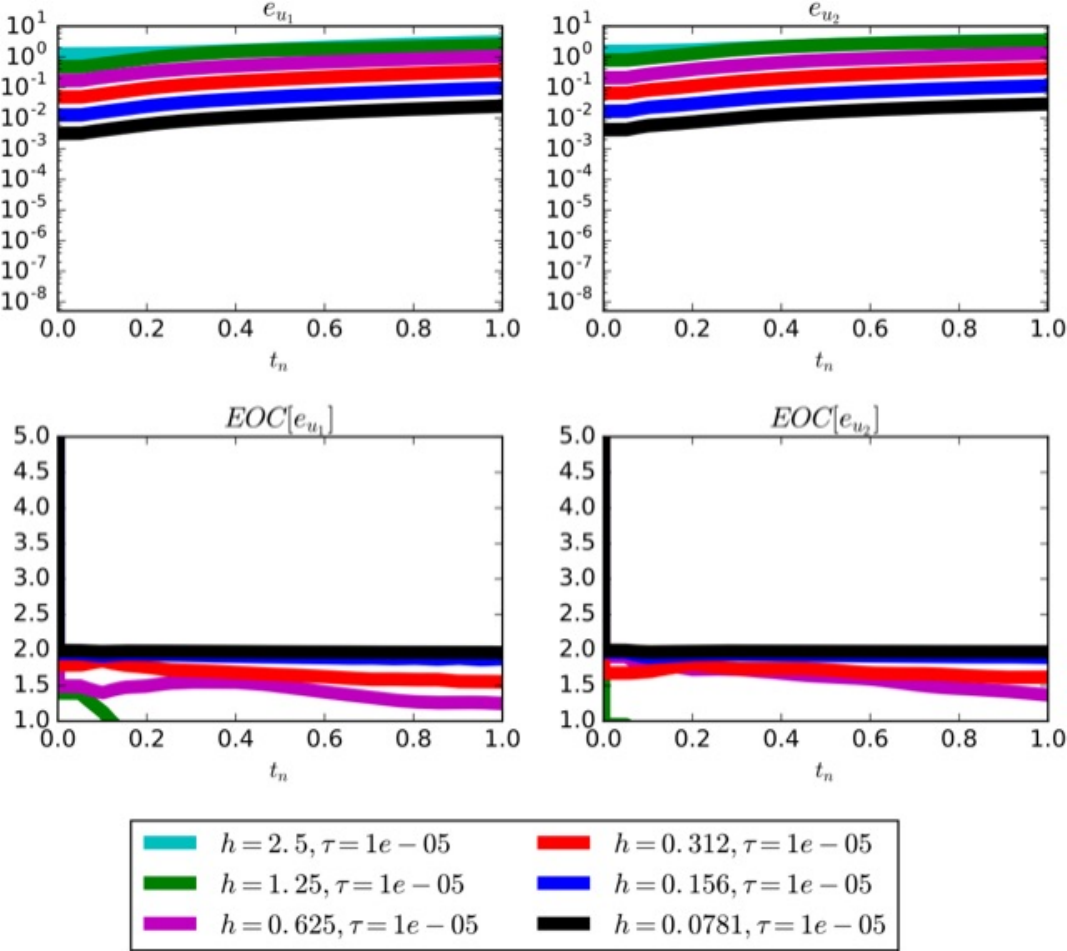}
 }
   \subfigure[][Here $q=2$ and we fix $\tau{} = 0.00001$. This is sufficiently small that the spatial discretisation error dominates.]{
   \includegraphics[scale=\figscale, width=0.47\figwidth]{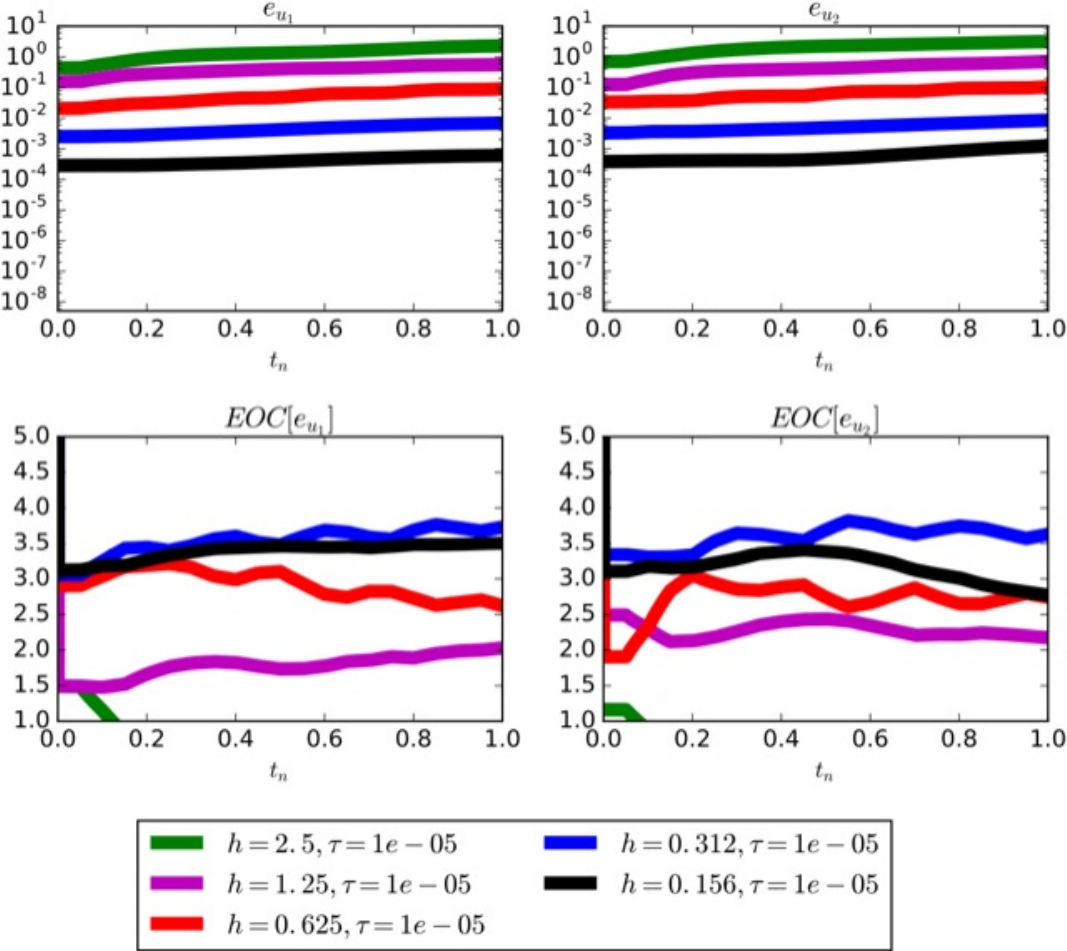}
 }
   \subfigure[][Here $q=3$ and we fix $\tau{} = 0.00001$. This is sufficiently small that the spatial discretisation error dominates.]{
   \includegraphics[scale=\figscale, width=0.47\figwidth]{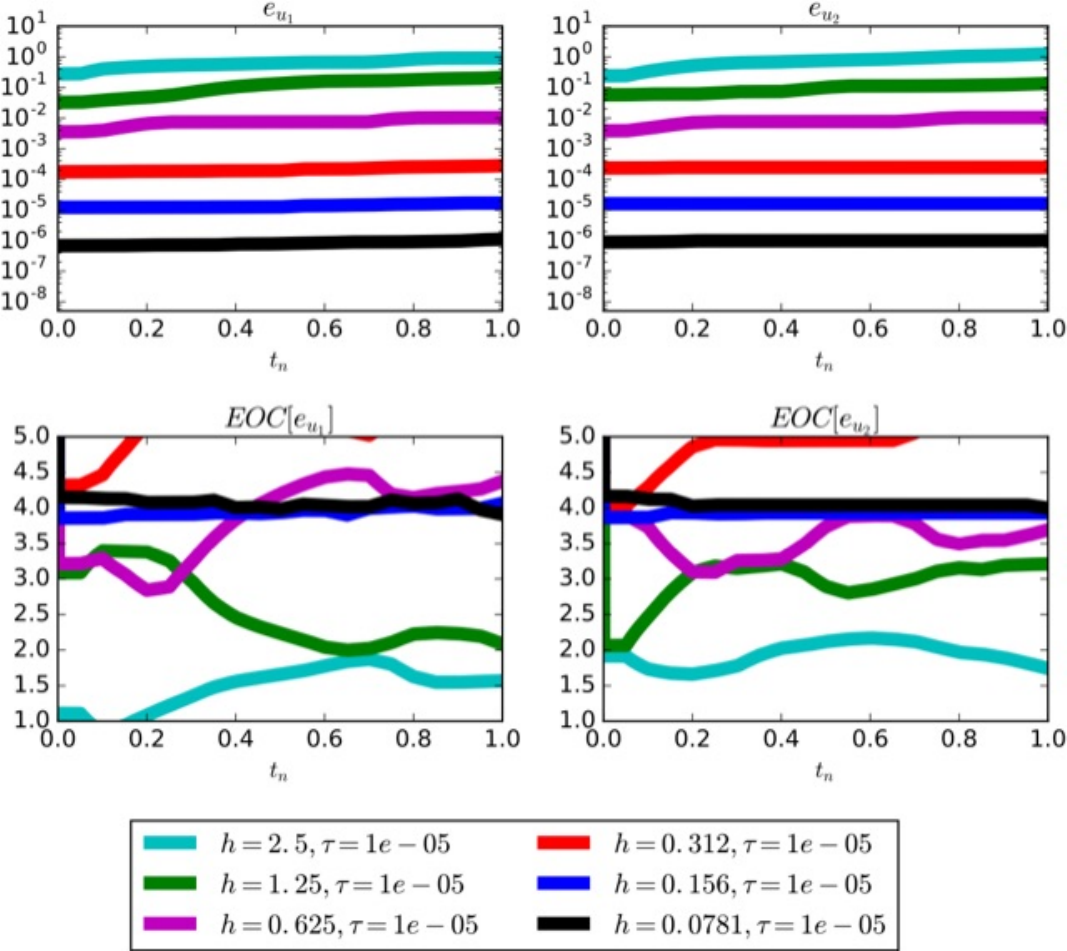}
 }
   \subfigure[][Here $q=2$ and on every refinement level we choose a coupling $\tau = C h$. Note that the time discretisation error here dominates.]{
     \includegraphics[scale=\figscale, width=0.47\figwidth]{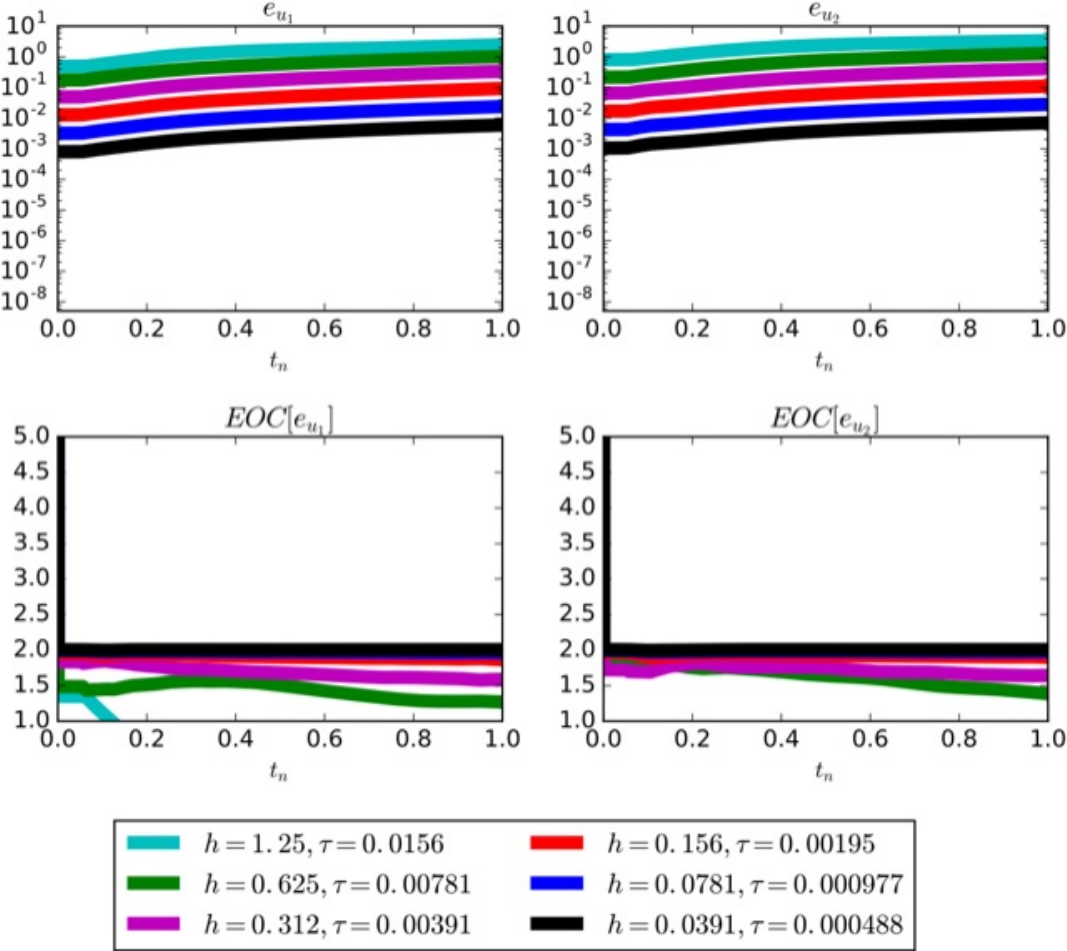}
   }
\end{figure}

We take $d=2$ and
\begin{equation}\label{eq:2solIC}
  \vec{u}_0=\frac{F_{\mu,\nu}}{G}\vec{E}_1+\frac{F_{\nu,\mu}}{G}\vec{E}_2,
\end{equation}
with $F_{\mu,\nu}$ given in (\ref{eq:F}) $G$ given in
(\ref{eq:G}). The parameters are $\vec{E_1} = (1,0)^\transpose$,
$\vec{E_2} = (0,1)^\transpose$, $\mu = \sqrt{2}$, $\nu = \sqrt{3}$,
$c_\nu = 24.9, c_\mu = 25.1$. The exact solution is then given by
(\ref{sec:exact-sol}). We take a uniform timestep and uniform meshes
that are fixed with respect to time. Convergence results are shown in
Figure \ref{fig:2solcon} and conservativity over long time is given in
Figure \ref{fig:2soldev}. Note that for $2$-soliton solution we have
$\vec W^{n} \neq \vec 0$ in general in which case the Lagrange
multiplier is required to ensure $\int_{\rS^1} \vec V^n \cdot \vec W^n
= 0$ for all $n$ and that the results of Theorem
\ref{the:cons-fullydiscrete} hold.

\clearpage
\subsection{Test 3 - Dynamics of $2$ and $3$-soliton solutions}

\subsubsection{Subtest 1}
\begin{figure}[h!]
 \centering
 \caption{
   \label{fig:2soldyn}
   Here we show the dynamics of the approximation generated by the conservative discretisation scheme with polynomial degree $q=1$ approximating a smooth solution with initial conditions given by (\ref{eq:2solICdyn}).}
   \subfigure{
   \includegraphics[scale=\figscale, width=0.31\figwidth]{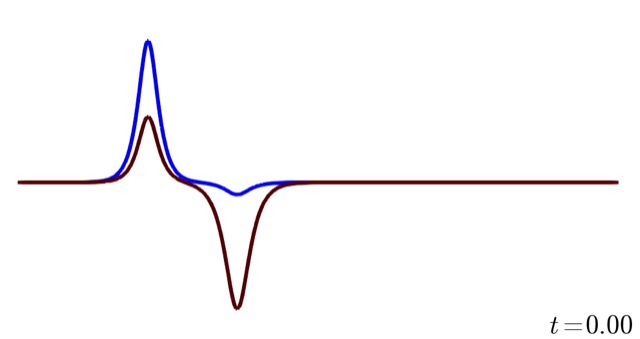}
 }
   \subfigure{
   \includegraphics[scale=\figscale, width=0.31\figwidth]{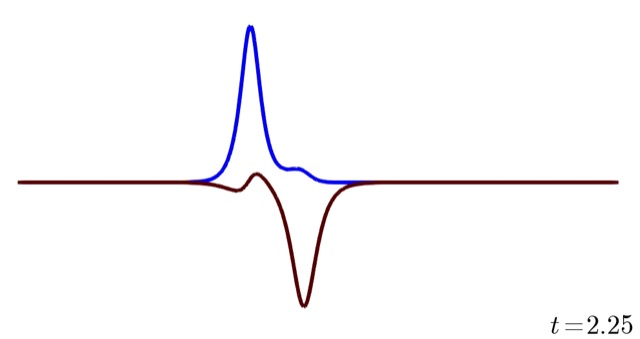}
 } 
  \subfigure{
   \includegraphics[scale=\figscale, width=0.31\figwidth]{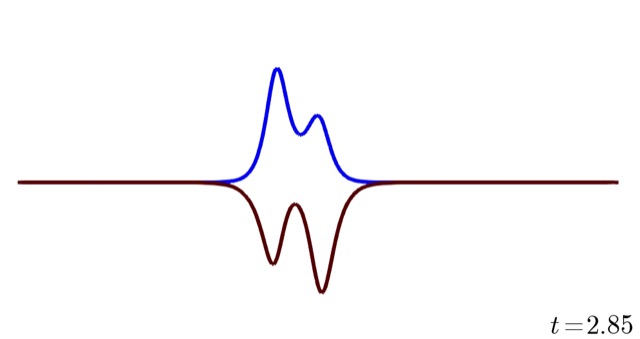}
  }
  \subfigure{
   \includegraphics[scale=\figscale, width=0.31\figwidth]{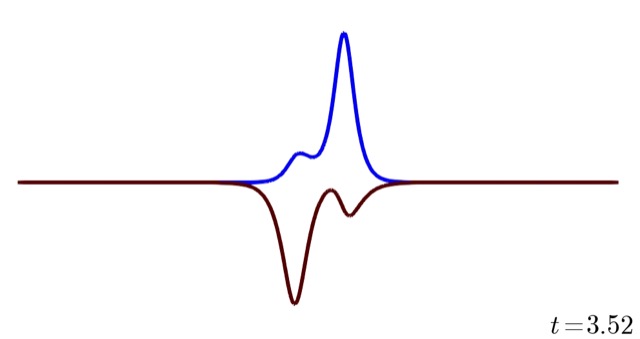}
 }
   \subfigure{
   \includegraphics[scale=\figscale, width=0.31\figwidth]{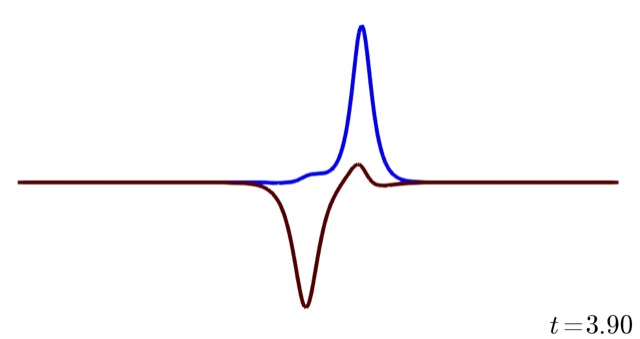}
 } 
  \subfigure{
   \includegraphics[scale=\figscale, width=0.31\figwidth]{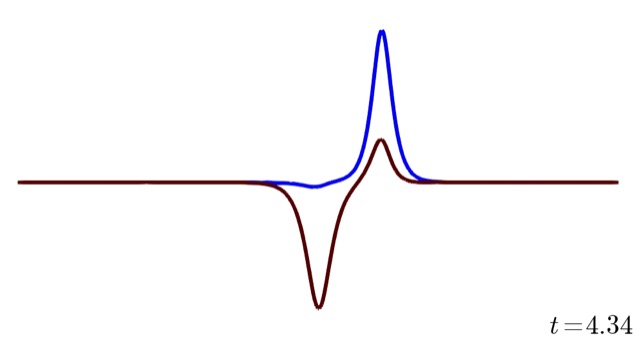}
 }
\end{figure}

We take $d=2$ and
\begin{equation}\label{eq:2solICdyn}
  \vec{u}_0=\frac{F_{\mu,\nu}}{G}\vec{E}_1+\frac{F_{\nu,\mu}}{G}\vec{E}_2,
\end{equation}
with $F_{\mu,\nu}$ given in (\ref{eq:F}) $G$ given in
(\ref{eq:G}). The parameters are $\vec{E_1} = (\frac 9
{10},\frac{\sqrt{19}}{10})^\transpose$, $\vec{E_2} = (\frac 1
{10},\frac{3\sqrt{11}}{10})^\transpose$, $\mu = \sqrt{2}$, $\nu =
\sqrt{3}$, $c_\nu = 10, c_\mu = 13$. Figure \ref{fig:2soldyn} shows
some plots of the dynamics of the numerical approximation.

\subsubsection{Subtest 2}
\begin{figure}[h!]
 \centering
 \caption{
   \label{fig:2soldyn2}
   Here we show the dynamics of the approximation generated by the conservative discretisation scheme with polynomial degree $q=1$ approximating a smooth solution with initial conditions given by (\ref{eq:2solICdyn2}).}
   \subfigure{
   \includegraphics[scale=\figscale, width=0.31\figwidth]{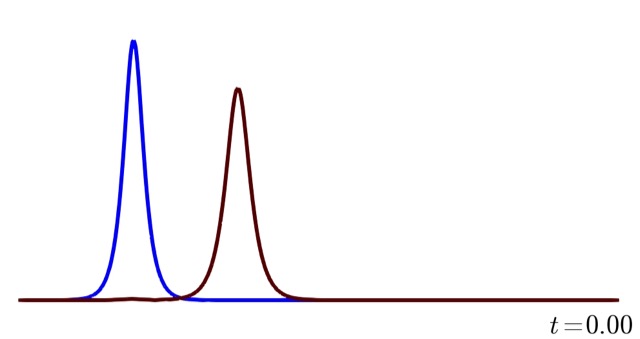}
 }
   \subfigure{
   \includegraphics[scale=\figscale, width=0.31\figwidth]{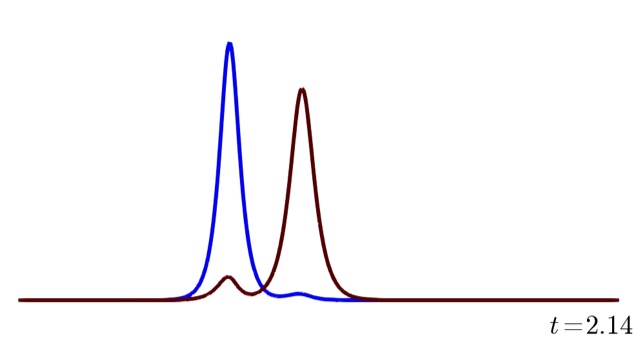}
 } 
  \subfigure{
   \includegraphics[scale=\figscale, width=0.31\figwidth]{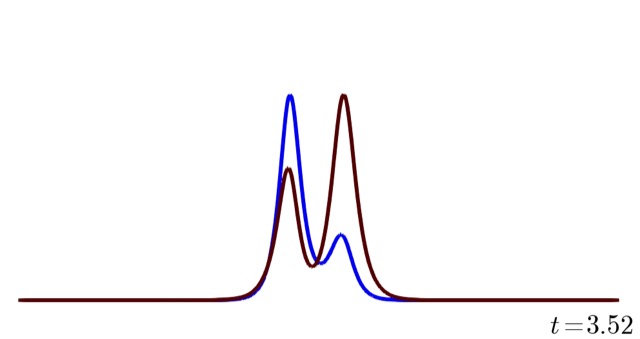}
  }
  \subfigure{
   \includegraphics[scale=\figscale, width=0.31\figwidth]{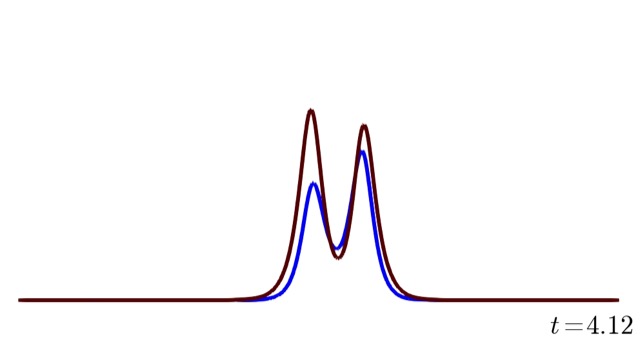}
 }
   \subfigure{
   \includegraphics[scale=\figscale, width=0.31\figwidth]{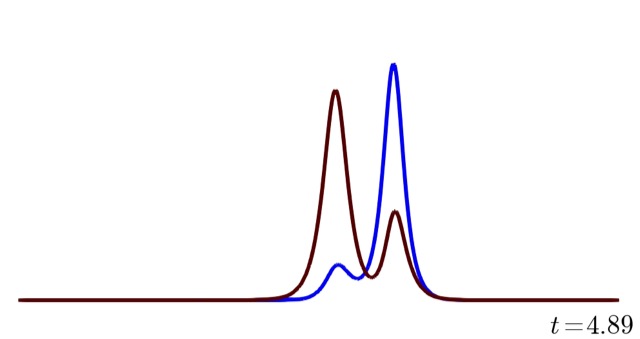}
 } 
  \subfigure{
   \includegraphics[scale=\figscale, width=0.31\figwidth]{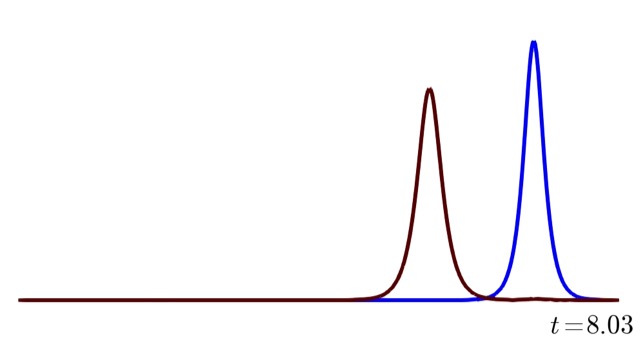}
 }
\end{figure}

We take $d=2$ and
\begin{equation}\label{eq:2solICdyn2}
  \vec{u}_0=\frac{F_{\mu,\nu}}{G}\vec{E}_1+\frac{F_{\nu,\mu}}{G}\vec{E}_2,
\end{equation}
with $F_{\mu,\nu}$ given in (\ref{eq:F}) $G$ given in
(\ref{eq:G}). The parameters are $\vec{E_1} = (1,0)^\transpose$,
$\vec{E_2} = (0,1)^\transpose$, $\mu =
\sqrt{2}$, $\nu = \sqrt{3}$, $c_\mu = 9, c_\nu = 13$. Figure
\ref{fig:2soldyn2} shows some plots of the dynamics of the numerical
approximation.

\subsubsection{Subtest 3}
\begin{figure}[h!]
 \centering
 \caption{
   \label{fig:3soldyn}
   Here we show the dynamics of the approximation generated by the conservative discretisation scheme with polynomial degree $q=1$ approximating a smooth solution with initial conditions given by (\ref{eq:3solICdyn2}).}
   \subfigure{
   \includegraphics[scale=\figscale, width=0.31\figwidth]{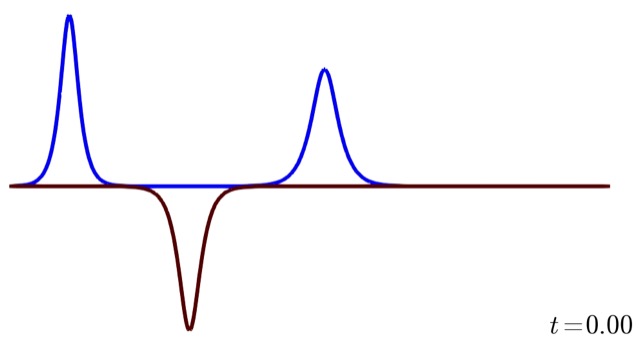}
 }
   \subfigure{
   \includegraphics[scale=\figscale, width=0.31\figwidth]{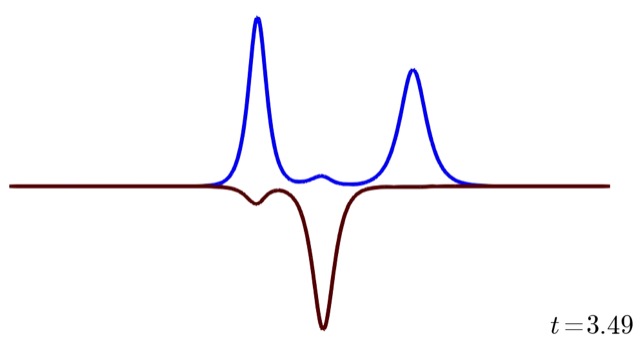}
 } 
  \subfigure{
   \includegraphics[scale=\figscale, width=0.31\figwidth]{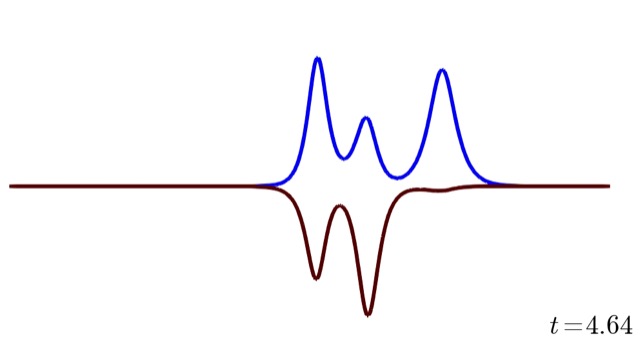}
  }
  \subfigure{
   \includegraphics[scale=\figscale, width=0.31\figwidth]{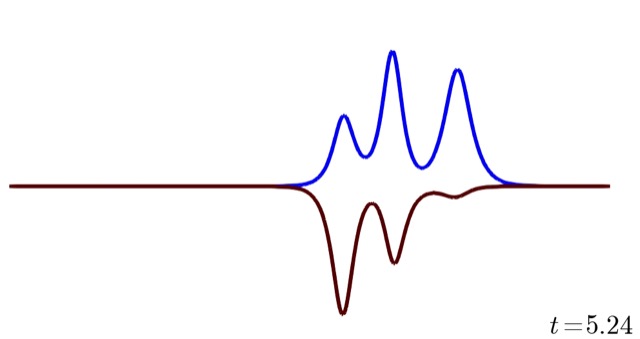}
 }
   \subfigure{
   \includegraphics[scale=\figscale, width=0.31\figwidth]{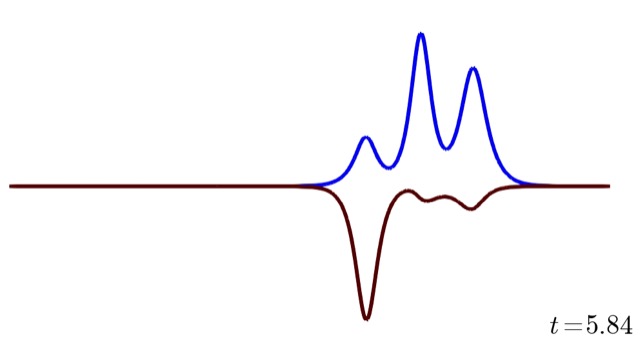}
 } 
  \subfigure{
   \includegraphics[scale=\figscale, width=0.31\figwidth]{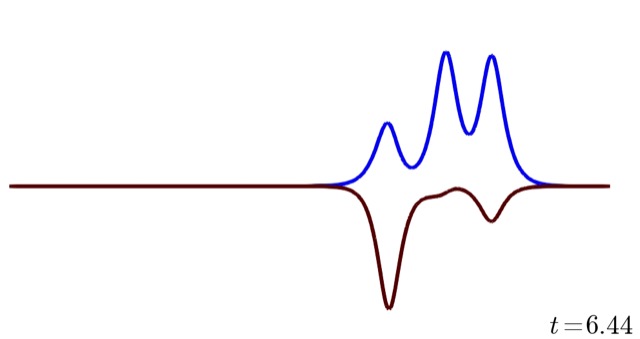}
 }
  \subfigure{
   \includegraphics[scale=\figscale, width=0.31\figwidth]{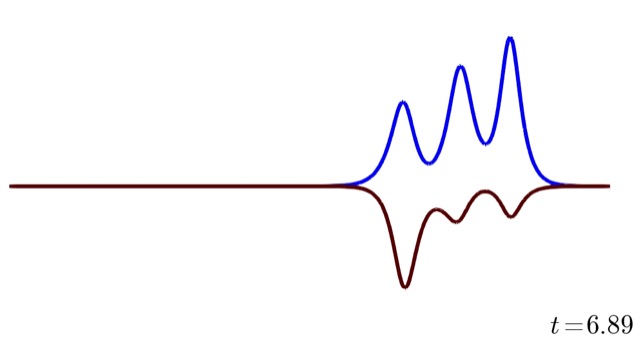}
 }
   \subfigure{
   \includegraphics[scale=\figscale, width=0.31\figwidth]{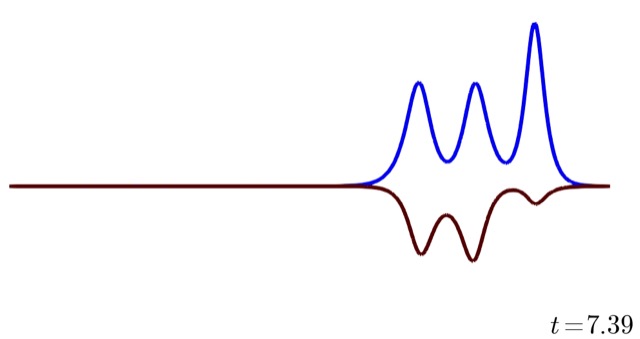}
 } 
  \subfigure{
   \includegraphics[scale=\figscale, width=0.31\figwidth]{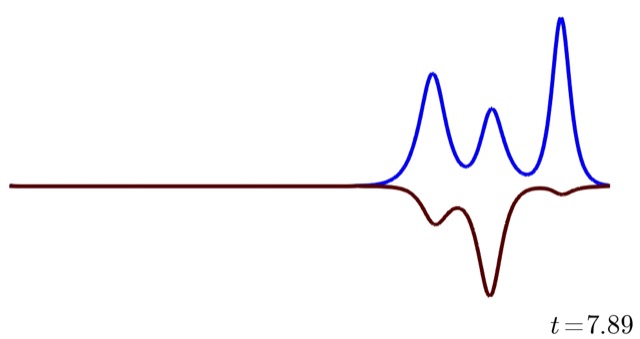}
 }
\end{figure}

In addition to the 2-soliton interactions we also take the opportunity
to examine the dynamics of a 3-soliton interaction. We take $d=2$ and
\begin{equation}
  \label{eq:3solICdyn2}
  \vec{u}_0
  =
  \sum_{i=1}^3
  \frac{2\mu_i}{\cosh\qp{\mu_i \qp{x - c_{\mu_i}}}}\vec{E_i}
\end{equation}
with $\vec E_1 = \vec E_3 = \qp{1,0}^\transpose, \vec E_2 = \qp{0,1}^\transpose$, $\mu_1 = \frac{19}{10}, \mu_2 = -\frac{40}{25}, \mu_3 = \frac{13}{10}$ and $c_{\mu_1} = 4, c_{\mu_2} = 12, c_{\mu_3} = 21$. Figure \ref{fig:3soldyn} shows the dynamics of the numerical approximation.

\clearpage

\subsection{Test 4 - Propagation of solitary waves from smooth initial data}

\begin{figure}[h]
 \centering
 \caption{
   \label{fig:trigcon}
   Here we examine the conservative discretisation scheme with various polynomial degrees, $q$, approximating the solution to (\ref{eq:vmkdv}) with initial conditions given by (\ref{eq:trig}). We show the deviation in the two invariants $F_i$, $i=2,4$, corresponding to momentum and energy respectively. In each test we take a fixed spatial discretisation parameter of $h=0.25$ and fixed time step of $\tau = 0.001$. Notice that in each case the deviation in energy is smaller than the solver tolerance of $10^{-12}$ and the deviation in momentum is bounded. The simulations are simulated for long time to test conservativity with $T=100$ in each case.}
   \subfigure[][$q=1$]{
   \includegraphics[scale=\figscale, width=0.31\figwidth]{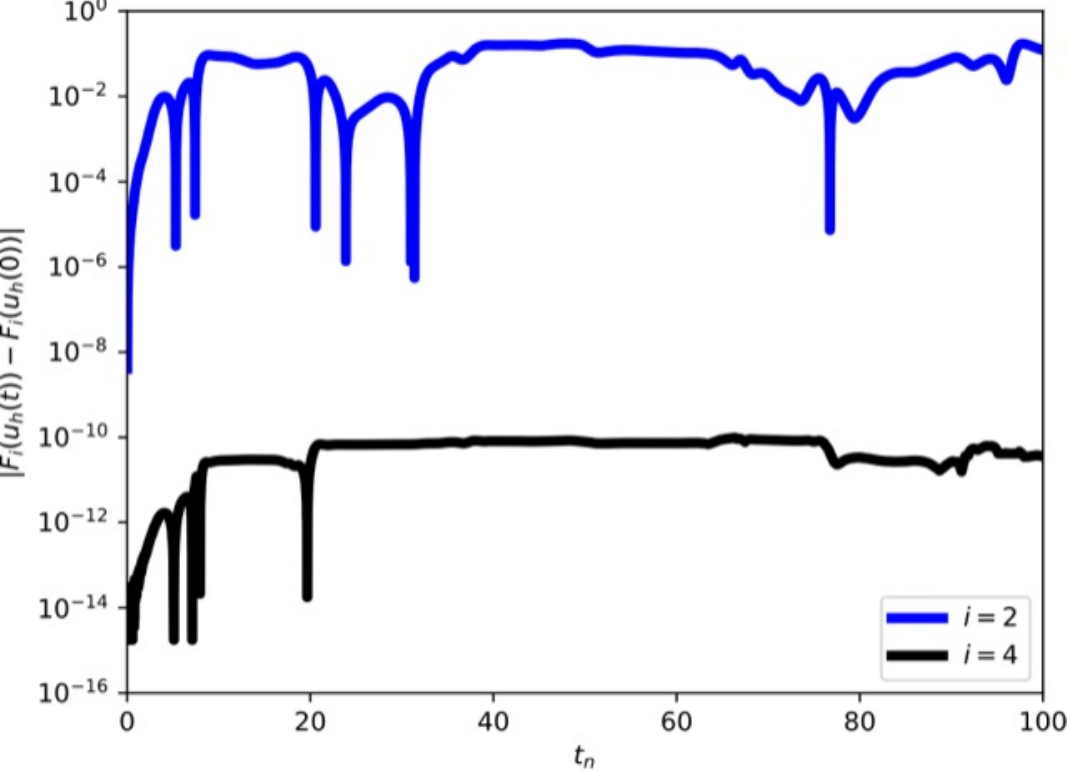}
 }
   \subfigure[][$q=2$]{
   \includegraphics[scale=\figscale, width=0.31\figwidth]{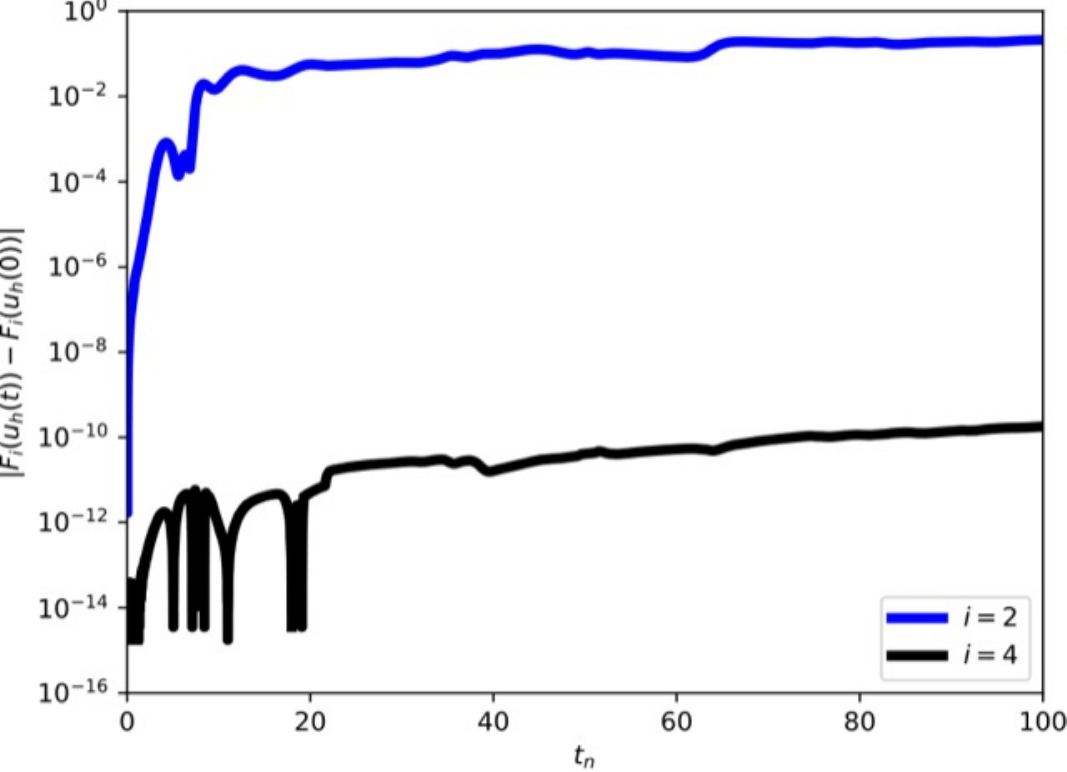}
 } 
  \subfigure[][$q=3$]{
   \includegraphics[scale=\figscale, width=0.31\figwidth]{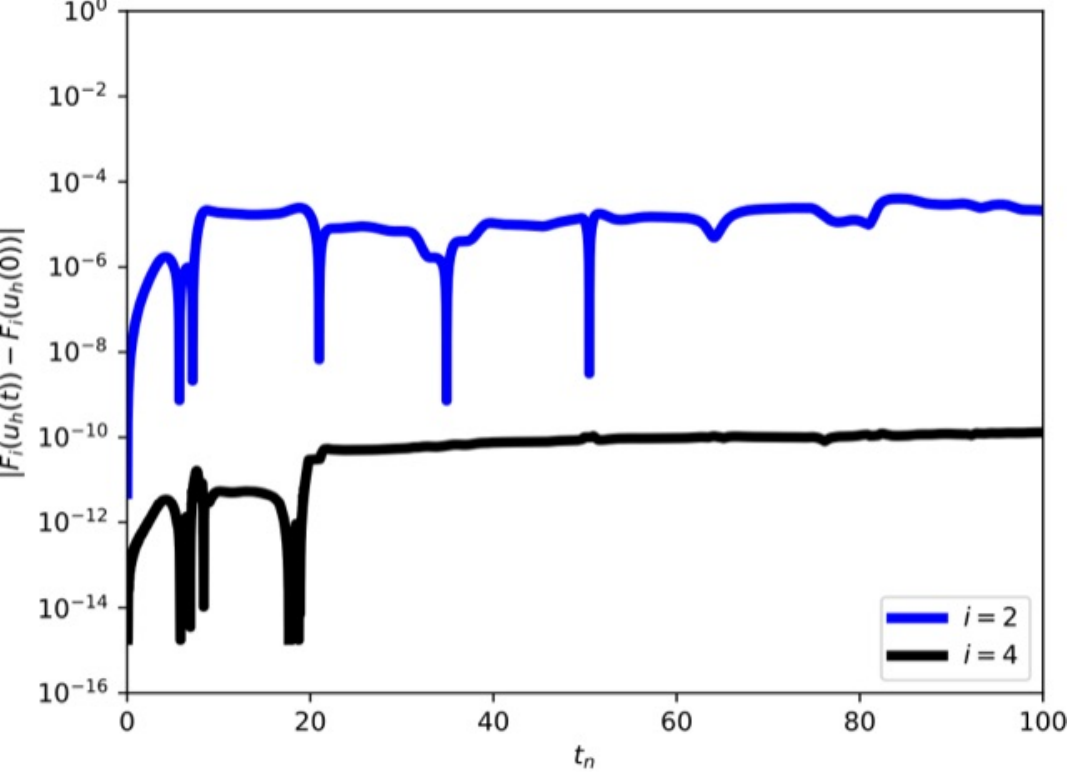}
 }
\end{figure}

\begin{figure}[h!]
 \centering
 \caption{
   \label{fig:trigdyn}
   Here we show the dynamics of the approximation generated by the conservative discretisation scheme with polynomial degree $q=1$ approximating the solution to (\ref{eq:vmkdv}) with initial conditions given by (\ref{eq:discon}). Notice that initially, dispersive waves emanate from the discontinuity.}
   \subfigure{
   \includegraphics[scale=\figscale, width=0.31\figwidth]{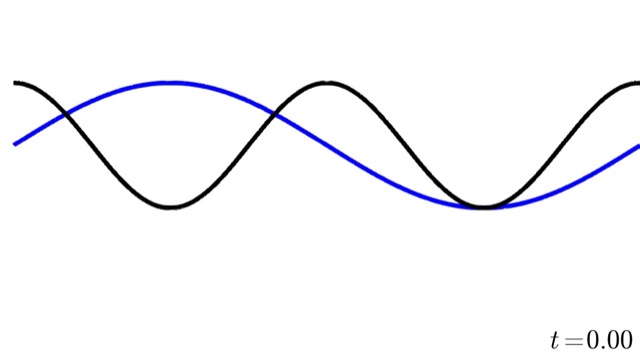}
 }
   \subfigure{
   \includegraphics[scale=\figscale, width=0.31\figwidth]{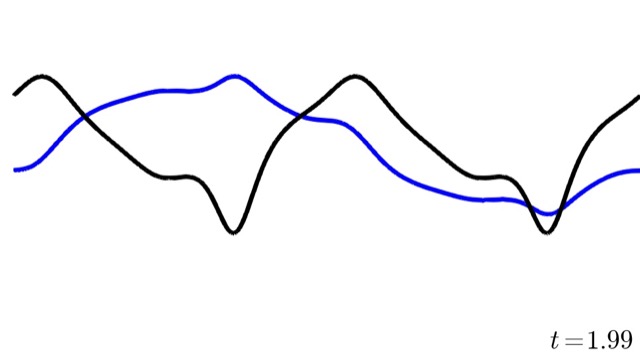}
 }
   \subfigure{
   \includegraphics[scale=\figscale, width=0.31\figwidth]{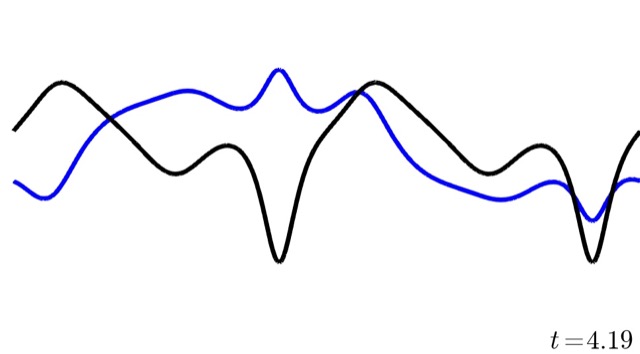}
 }
   \subfigure{
   \includegraphics[scale=\figscale, width=0.31\figwidth]{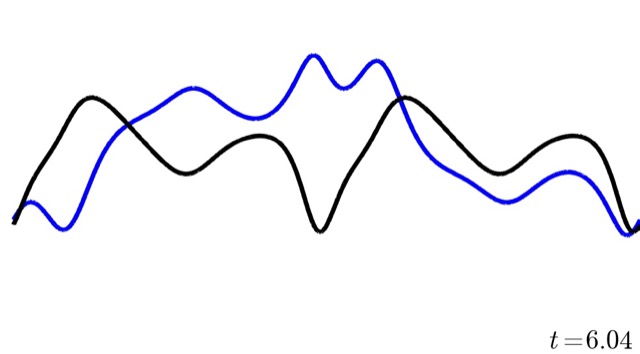}
 }
   \subfigure{
   \includegraphics[scale=\figscale, width=0.31\figwidth]{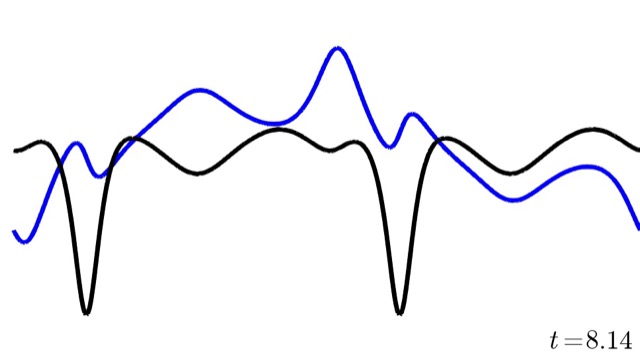}
 }
   \subfigure{
   \includegraphics[scale=\figscale, width=0.31\figwidth]{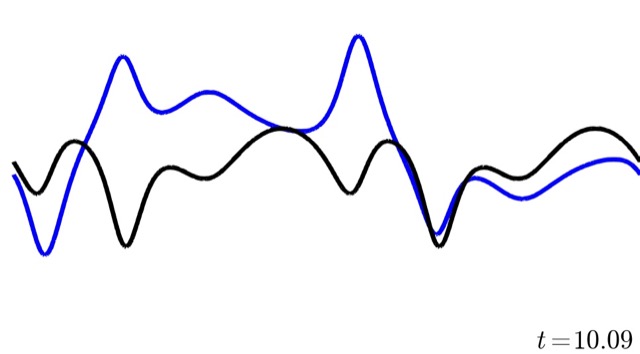}
 }
   \subfigure{
   \includegraphics[scale=\figscale, width=0.31\figwidth]{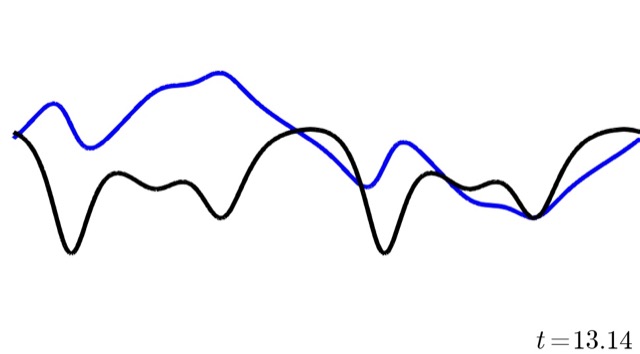}
 }
   \subfigure{
   \includegraphics[scale=\figscale, width=0.31\figwidth]{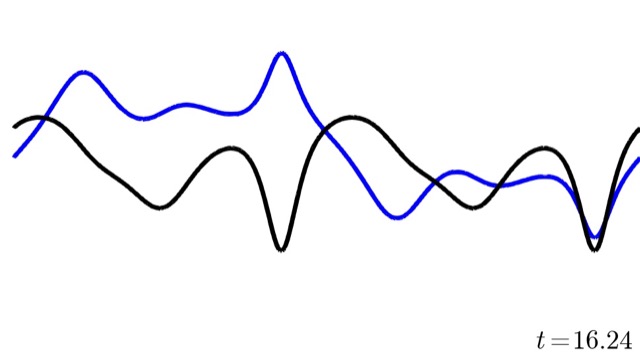}
 }
   \subfigure{
   \includegraphics[scale=\figscale, width=0.31\figwidth]{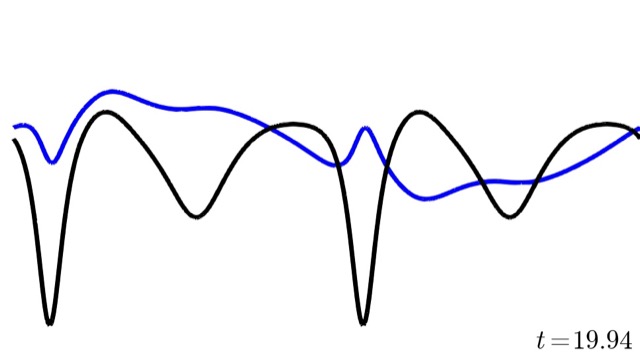}
 }
\end{figure}

We take $d=2$ and $\vec u_0 = \qp{u_{0,1}, u_{0,2}}$ with
\begin{equation}
  \label{eq:trig}
  \begin{split}
    u_{0,1}
    &=
    \sin{\frac{\pi}{20}x}
    \\
    u_{0,2}
    &=
    \cos{\frac{\pi}{10}x}.
  \end{split}
\end{equation}
The solution here is smooth and solitary waves begin to form quickly
into the simulation. Plots of the solutions are given in Figure
\ref{fig:trigdyn} as well as conservativity plots in Figure
\ref{fig:trigcon}.

\clearpage
\subsection{Test 5 - Solution with discontinuous initial data}

\begin{figure}[h]
 \centering
 \caption{
   \label{fig:discon}
   Here we examine the conservative discretisation scheme with various polynomial degrees, $q$, approximating the solution to (\ref{eq:vmkdv}) with initial conditions given by (\ref{eq:discon}). We show the deviation in the two invariants  $F_i$, $i=2,4$, corresponding to momentum and energy respectively. In each test we take a fixed spatial discretisation parameter of $h=0.25$ and fixed time step of $\tau = 0.001$. Notice that in each case the deviation in energy is smaller than the solver tolerance of $10^{-12}$ and the deviation in momentum is bounded. The simulations are simulated for long time to test conservativity with $T=100$ in each case.}
   \subfigure[][$q=1$]{
   \includegraphics[scale=\figscale, width=0.31\figwidth]{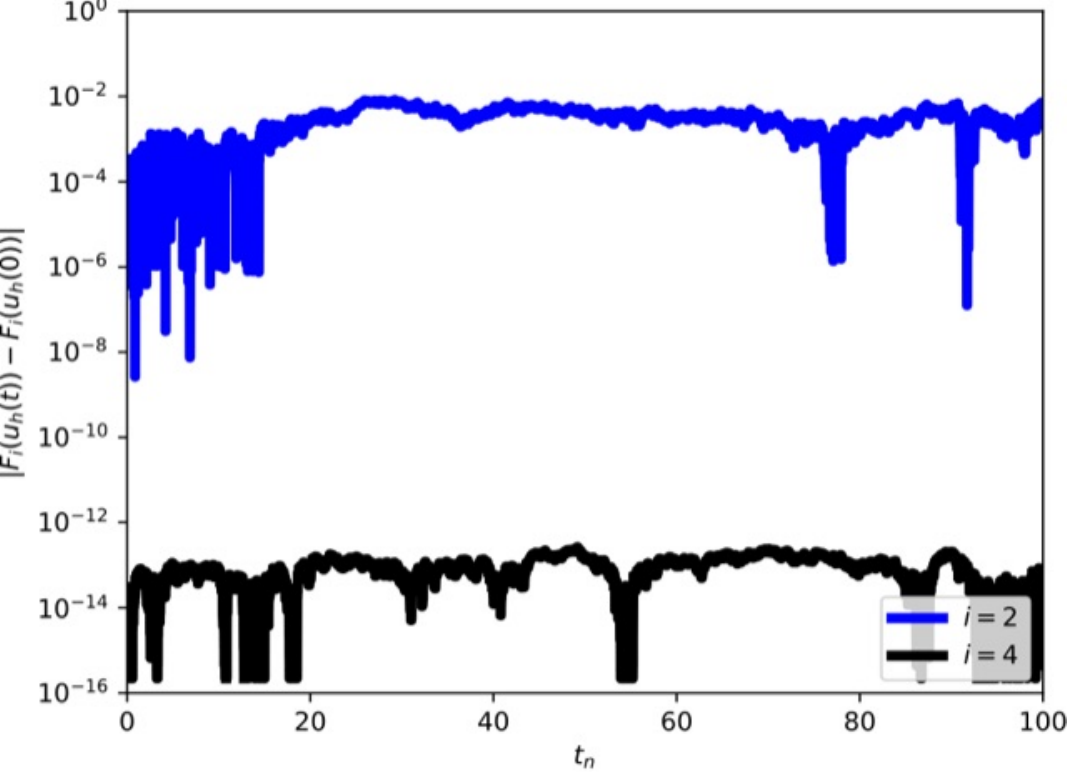}
 }
   \subfigure[][$q=2$]{
   \includegraphics[scale=\figscale, width=0.31\figwidth]{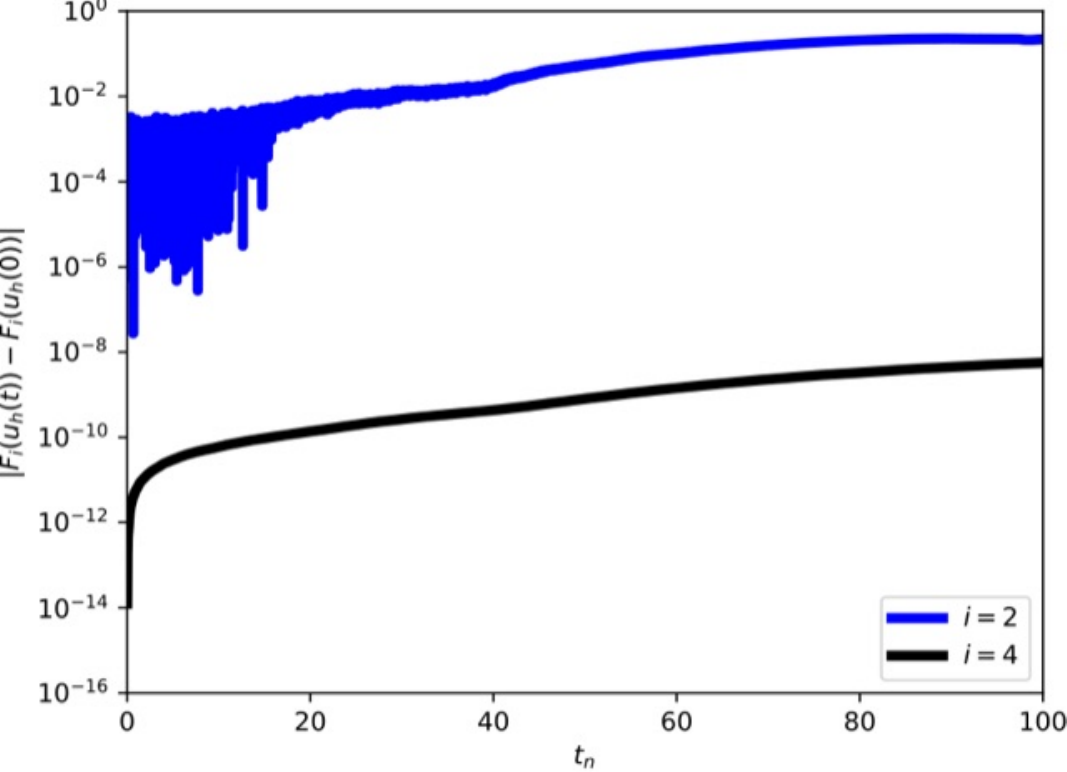}
 } 
  \subfigure[][$q=3$]{
   \includegraphics[scale=\figscale, width=0.31\figwidth]{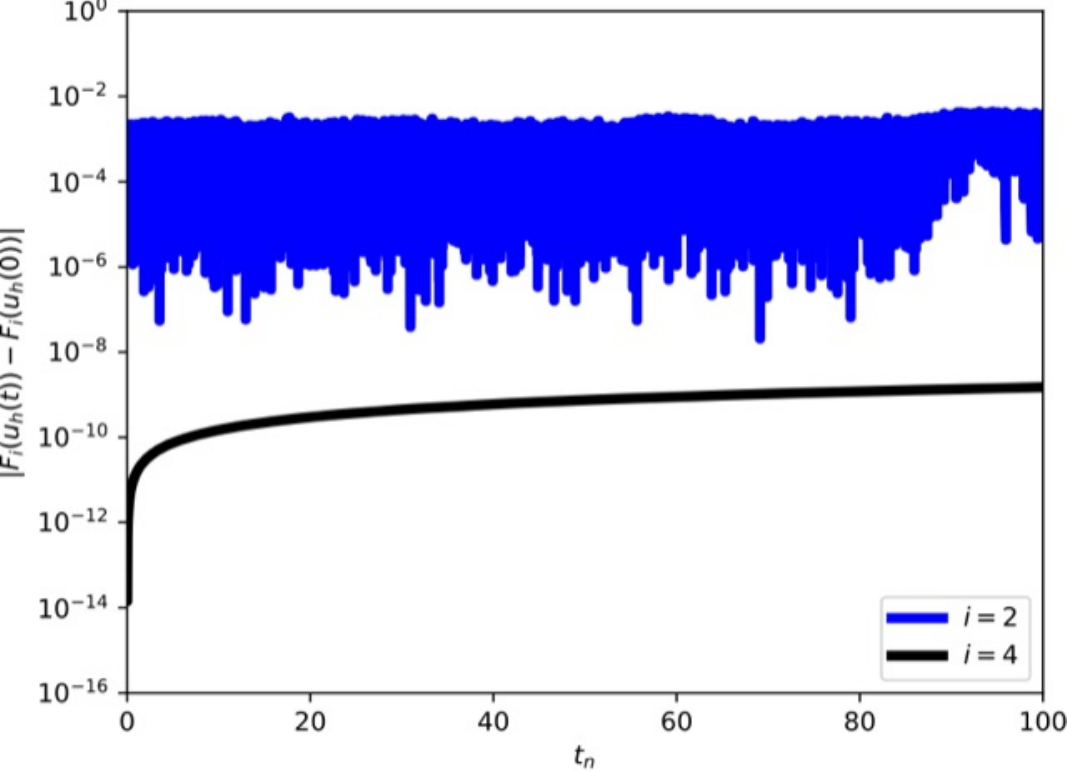}
 }
\end{figure}

\begin{figure}[h!]
 \centering
 \caption{
   \label{fig:disdyn}
   Here we show the dynamics of the approximation generated by the conservative discretisation scheme with polynomial degree $q=1$ approximating the solution to (\ref{eq:vmkdv}) with initial conditions given by (\ref{eq:discon}). Notice that initially, dispersive waves emanate from the discontinuity.}
   \subfigure{
   \includegraphics[scale=\figscale, width=0.31\figwidth]{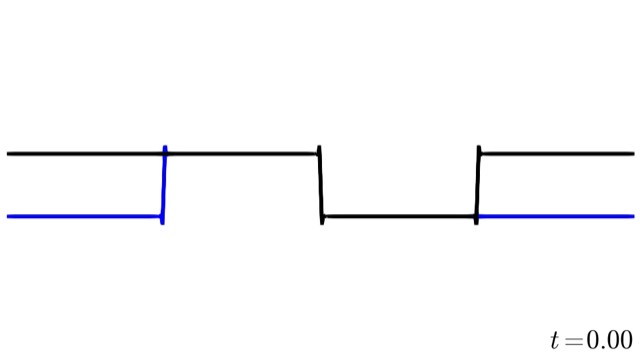}
 }
   \subfigure{
   \includegraphics[scale=\figscale, width=0.31\figwidth]{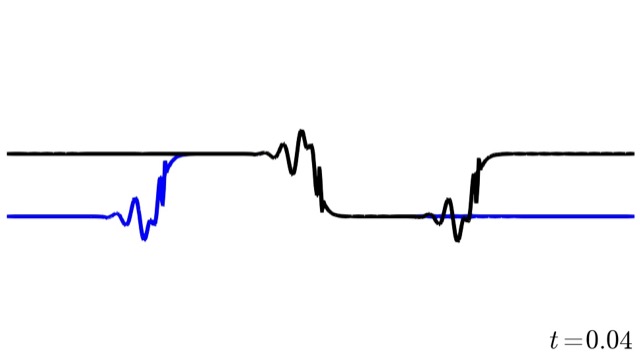}
 }
   \subfigure{
   \includegraphics[scale=\figscale, width=0.31\figwidth]{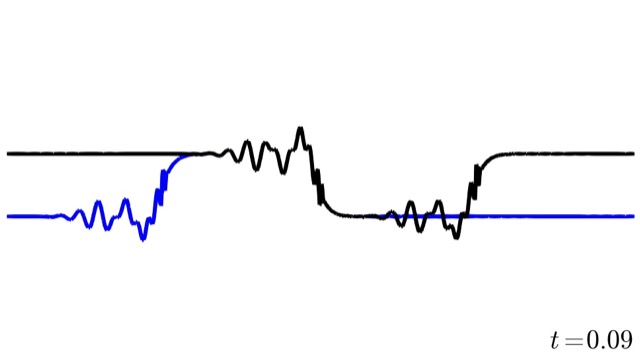}
 }
   \subfigure{
   \includegraphics[scale=\figscale, width=0.31\figwidth]{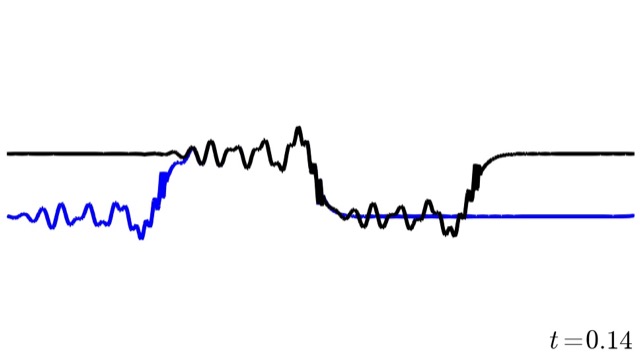}
 }
   \subfigure{
   \includegraphics[scale=\figscale, width=0.31\figwidth]{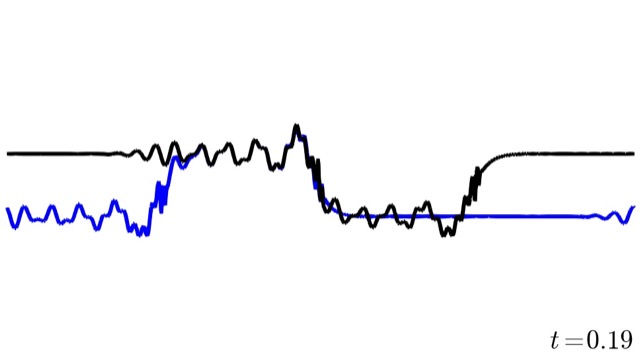}
 }
   \subfigure{
   \includegraphics[scale=\figscale, width=0.31\figwidth]{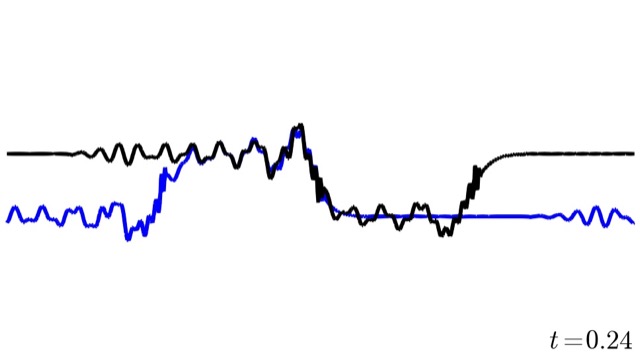}
 }
   \subfigure{
   \includegraphics[scale=\figscale, width=0.31\figwidth]{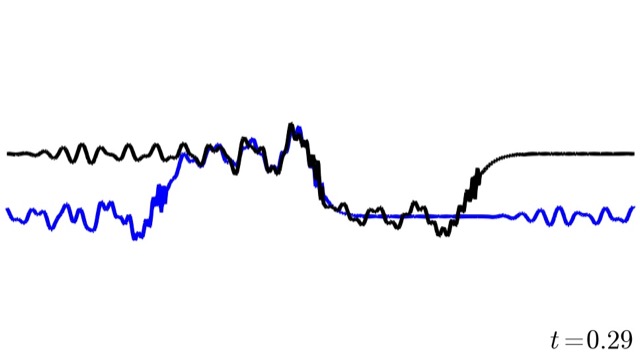}
 }
   \subfigure{
   \includegraphics[scale=\figscale, width=0.31\figwidth]{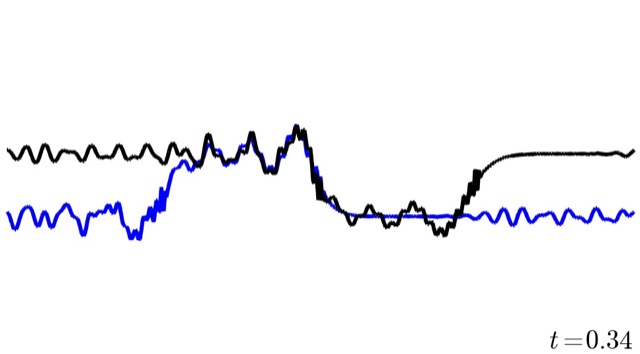}
 }
   \subfigure{
   \includegraphics[scale=\figscale, width=0.31\figwidth]{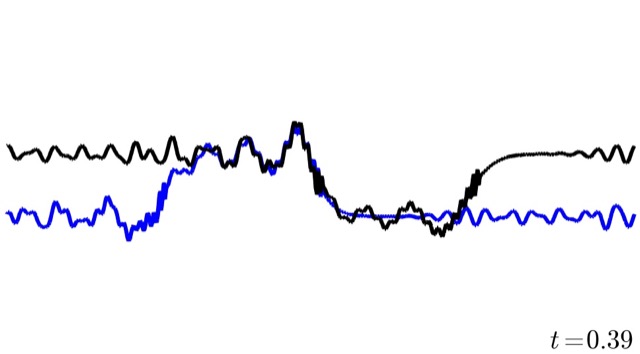}
 }
   \subfigure{
   \includegraphics[scale=\figscale, width=0.31\figwidth]{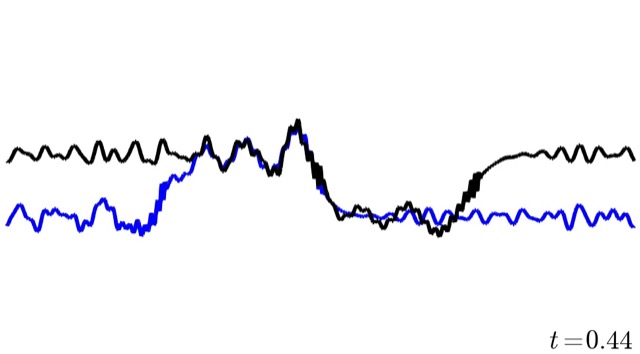}
 }
   \subfigure{
   \includegraphics[scale=\figscale, width=0.31\figwidth]{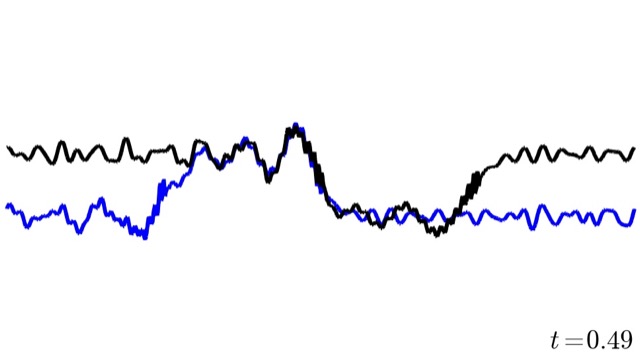}
 }
   \subfigure{
   \includegraphics[scale=\figscale, width=0.31\figwidth]{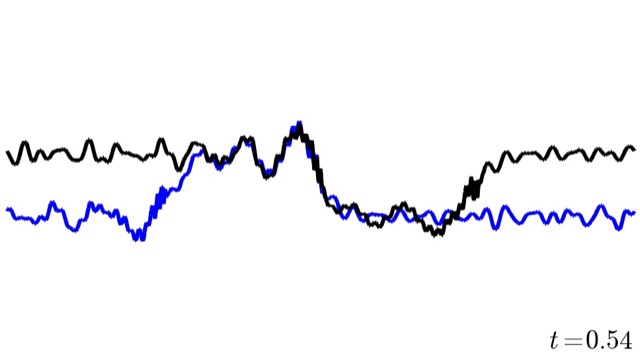}
 }
\end{figure}

We take $d=2$ and $\vec u_0 = \qp{u_{0,1}, u_{0,2}}$ with
\begin{equation}
  \label{eq:discon}
  \begin{split}
    u_{0,1}
    &=
    \begin{cases}
      1 \text{ for } x\in [10, 20]
      \\
      0 \text{ otherwise}.
    \end{cases}
    \\
    u_{0,2}
    &=
    \begin{cases}
      0 \text{ for } x\in [20, 30]
      \\
      1 \text{ otherwise}.
    \end{cases}
  \end{split}
\end{equation}
The solution here is discontinuous in both components. This is a
particularly tough scenario to simulate as there is no guarantee of
classical solutions. We align the mesh to the discontinuities so that
the discrete energy at the initial condition makes sense. Plots of the
solutions are given in Figure \ref{fig:disdyn} as well as conservativity plots
in Figure \ref{fig:discon}.

\section{Conclusions}

In this work we have constructed a consistent Galerkin approximation
for the \vmkdv equation. We have proven that both the
semi-discretisations as well as the fully discrete problems are
conservative and numerically shown that this is true in practice. In
addition, we have given numerical evidence to suggest that the method
is of optimal order, that is, $\Norm{\vec u - \vec
  U}_{\leb{\infty}((0,T); \leb{2}(\rS^1))} = \Oh(\tau^2 +
h^{q+1})$. We expect methods designed in this fashion, which is quite
generic, to be successful in the simulation of geophysical fluid flows.

\clearpage
\bibliographystyle{alpha}
\bibliography{./tristansbib,./tristanswritings,./central}

\end{document}